\documentclass[11pt,a4paper,leqno]{amsart} 
 
\usepackage{amsmath,amssymb,amsfonts,amsthm,graphics} 
\usepackage[matrix,arrow,curve]{xy}

\newcommand{\DP}[2]{{\frac{\partial #1}{\partial #2}}}  
 
\newcommand{\tx}[1]{\text{\rm #1}} 
     
\newcommand{\CC}{\mathbb{C}}    
   
\newcommand{\QQ}{\mathbb{Q}}    
\newcommand{\RR}{\mathbb{R}}    
\newcommand{\Ss}{\mathbb{S}}

\newcommand{\ZZ}{\mathbb{Z}}

\newcommand{\demo}{\noindent {\it \small Proof. \ }}    
  
\newtheorem{defi}{\sf Definition}[section]   
\newtheorem{thm}[defi]{\sf Theorem}      
\newtheorem{cor}[defi]{\sf Corollary}    
    
\newtheorem{lem}[defi]{\sf Lemma}

\def\mc{\mathcal}

\title[Floer homology for manifolds with  boundary or
Symplectic homology]
{A survey of Floer homology for manifolds with contact type boundary 
or Symplectic homology}  
\author{Alexandru OANCEA} 
\date{March 4, 2004 - Revised version.}

\begin{document} 
 
\maketitle

\begin{abstract} 

  The purpose of this paper is to give a survey of the various
  versions of Floer homology for manifolds with contact type boundary 
  that have so far appeared  in the
  literature. Under the name of 
  ``Symplectic homology'' or ``Floer homology for manifolds with
  boundary'' they bear in fact common  features and we shall try to
  underline the  principles that unite them. Once this will be
  accomplished we shall proceed to describe the 
  peculiarity of each of the constructions and the specific
  applications that unfold out of it: classification of ellipsoids and
  polydiscs in $\CC^n$, stability of the action spectrum for 
  contact type boundaries of symplectic manifolds, existence of closed
  characteristics on contact type hypersurfaces and obstructions to
  exact Lagrange embeddings. The computation of
  the Floer cohomology for balls in $\CC^n$ is carried by explicitly
  perturbing the nondegenerate Morse-Bott spheres of closed
  characteristics. 
\end{abstract}

\tableofcontents

\renewcommand{\thefootnote}{\fnsymbol{footnote}}
\setcounter{footnote}{0}
\footnotetext{ 
{\it Keywords}: Floer homology and cohomology -
Hamiltonian systems.

{\it 2000 Mathematics Subject Classification}: 53D40 - 37J45.

Work supported by ENS Lyon, Ecole Polytechnique (France) and ETH
Zurich (Switzerland).}   
\renewcommand{\thefootnote}{\arabic{footnote}}
\setcounter{footnote}{0}

\renewcommand{\labelitemi}{$-$}

\newcounter{firstTSW} 
\addtocounter{equation}{1} 
\setcounter{firstTSW}{\value{equation}}
 
\newcounter{secondTSW} 
\addtocounter{equation}{1} 
\setcounter{secondTSW}{\value{equation}}
 
\newcounter{thirdTSW} 
\addtocounter{equation}{1} 
\setcounter{thirdTSW}{\value{equation}}

\section{Constructions of Floer homological invariants} 
 
\subsection{Morse homology} \label{subsection 1} 

Floer homology can be seen as a vast
generalization of the Thom-Smale-Witten (TSW) complex to infinite
dimension. We briefly recall below the construction of the
latter \cite{Fl, Sch}. The resulting homology theory is called Morse 
(co)homology as a
tribute to Marston Morse's pioneering use of critical points of smooth
maps in the study of the topology of manifolds.

Consider a Morse function $f: L \longrightarrow \RR$ defined
on a closed manifold $L$ and a Riemannian metric $g$ which is generic
with respect to a certain transversality property. One defines two
complexes, one homological and the other one cohomological, by

\bigskip

{\scriptsize 
\begin{tabular}{@{}p{220pt}|p{220pt}} 
(\thesection.\thefirstTSW) \  
$\displaystyle  
C_q^{\tx{Morse}}(L,-\nabla^g  f) = \bigoplus _{\scriptsize \begin{array}{c} 
    \nabla^g  f (x) =0 \\ \tx{ind}_f(x)=q \end{array} }  
\ZZ \langle x \rangle \, , $  
&  
\ (\thesection.\thefirstTSW$^\prime$) \  
$\displaystyle  
C^q_{\tx{Morse}}(L,-\nabla^g  f) = \bigoplus _{\scriptsize \begin{array}{c} 
    \nabla^g  f (x) =0 \\ \tx{ind}_f(x)=q \end{array} }  
\ZZ \langle x 
\rangle \, , $ \\ 
 & \\ & \\  
(\thesection.\thesecondTSW) \ 
$\partial ^{\tx{Morse}}~: C_q^{\tx{Morse}}(L,-\nabla^g  f) 
\longrightarrow C_{q-1}  
^{\tx{Morse}}(L,-\nabla^g  f) \, ,$ 
& 
\ (\thesection.\thesecondTSW$^\prime$) \  
$\delta _{\tx{Morse}}~: C^q_{\tx{Morse}}(L,-\nabla^g  f) 
\longrightarrow C^{q+1}  
_{\tx{Morse}}(L,-\nabla^g  f) \, ,$ \\ 
 & \\ & \\ 
(\thesection.\thethirdTSW) \ 
$\displaystyle  
\partial ^{\tx{Morse}}\langle x \rangle = \sum _{\tx{ind}_f(y)=q-1} 
\# \big( \mathcal {M}(x,\, y~; \, -\nabla^g  f)  
/ \RR \big) \ \langle y \rangle \, . $ 
& 
\ (\thesection.\thethirdTSW$^\prime$) \  
$\displaystyle  
\delta _{\tx{Morse}}\langle x \rangle = \sum _{\tx{ind}_f(y)=q+1} 
\# \big( \mathcal {M}(y,\, x~; \, -\nabla^g  f)  
/ \RR \big) \ \langle y \rangle \, . $ 
\end{tabular}  
} 
 
\bigskip  
 
The notation $\mathcal {M}(y,\, x~; \, -\nabla^g  f) $ stands for the
space of parameterized negative gradient trajectories running from
$y$ to $x$, on which the additive group $\RR$ naturally acts by
reparameterization. 
The transversality condition mentioned above concerns the transverse
intersection of any pair of stable and unstable manifolds of $-\nabla^g  
f$. It ensures in particular that  $\mathcal {M}(y,\, x~; \, 
-\nabla^g  f) $ is a smooth manifold of dimension 
$\tx{ind}_f(y) - \tx{ind}_f(x)$, where $\tx{ind}_f(x)$ 
is the  Morse index of the nondegenerate 
critical point $x$ for the function $f$. 
When the difference of the indices is equal to $1$, a
careful description of the relative compactness of sequences of
trajectories proves that the 
quotient  $\mathcal{M}(y,\, x~; \, -\nabla^g  f)  / \RR  $ is finite.  
Any choice of orientations for the unstable manifolds allows one to
orient the trajectory spaces and algebraically count the elements of 
$\mathcal{M}(y,\, x~; \, -\nabla^g  f)  / \RR  $ when $\tx{ind}_f(y) -
\tx{ind}_f(x)$ is equal to $1$. This gives a precise meaning for the
expressions defining 
$\partial ^{\tx{Morse}}$ and  
$\delta_{\tx{Morse} } $. The crucial identities  
$\big( \partial^{\tx{Morse}} \big) ^2=0$ and  
$\big( \delta  _{\tx{Morse}} \big) ^2 =0$ 
are a consequence of a glueing theorem
which constitutes, together with the analysis of the convergence of
sequences of trajectories, a description of the compactification of 
$\mathcal{M}(y,\, x; \, -\nabla^g  f)/ \RR$ by ``broken trajectories''.  
Finally, the homology of the above complexes is seen to compute the
singular homology and, respectively, cohomology 
$H_*(C_*^{\tx{Morse}}(L, \, -\nabla^g  f) ) 
\simeq H_*(L~; \, \ZZ)$, $H^*(C^*_{\tx{Morse}}(L, \, -\nabla^g  f))  
\simeq H^*(L~; \, \ZZ)$. The most intuitive way to see this is to use
the fact \cite{La} that the unstable manifolds give rise to a
CW-decomposition of $L$ and use cellular (co)homology as a bridge
between Morse and singular (co)homology \cite{Po, teza mea}. 

We should mention at this point two important extensions of the
preceding setup. Firstly, the negative gradient vector field can be
replaced by a vector field $X$ which is {\it negative gradient-like} 
with
respect to a Morse function 
$f$ i.e. $X\cdot f <0$ outside the critical points of $f$,
 which satisfies the same transversality condition and which 
is equal to
$-\nabla ^g f$ near the critical points for some given metric $g$. 
Secondly, when the
manifold $L$ is no longer compact, the same construction can be
carried out if $f$ satisfies the {\it Palais-Smale} condition: every
sequence $(x_\nu)$ such that $f(x_\nu)$ is bounded and $\nabla^g  f(x_\nu) 
\longrightarrow 0$  contains a convergent subsequence. The
Palais-Smale condition can of course be formulated for a (negative)
gradient-like vector field. For any two regular values $a<b$ of $f$ 
the complexes  
constructed on the critical points belonging to  $f^ 
{-1} [a,\, b]$ eventually compute the relative singular
homology/cohomology 
$H_*(C_*^{\tx{Morse}}(L,\, -\nabla^g  f~; \, a, \, b ) ) 
\simeq H_*(f^b, \, f^a;\, \ZZ)$,   
$H^*(C^*_{\tx{Morse}}(L,\, -\nabla^g  f~; \, a, \, b ) ) 
\simeq H^*(f^b, \, f^a;\, \ZZ)$. For any real number $c$ we have
denoted  $f^{c} = \{ f \le c \}$.

\bigskip 
 
\subsection{Floer homology for closed manifolds} 
\label{hom Floer var fermees} 
Let us now describe the main lines along which Hamiltonian Floer
homology of a closed manifold is constructed. This section recalls the
construction presented in F. Laudenbach's paper ``Symplectic geometry
and Floer homology'' published in this same issue of the journal. 
A comprehensive reference is 
also provided by D. Salamon's lecture notes \cite{Sa}, while the full
details  under the symplectic asphericity condition 
$\langle \omega , \, \pi_2(M) \rangle = 0$ that we require henceforth are
disseminated in several papers of A. Floer
\cite{F1, F2}. This assumption can be eliminated with moderate effort
in the monotone case (``bubbles'') and with strenuous one (``virtual
cycle technique'') in the general case, but it is suitable for this
survey paper to simplify things. 
We henceforth assume $\langle \omega, \, \pi_2(M) \rangle =0$, as well as  
$\langle c_1, \, \pi_2(M) \rangle =0$ in order for 
 the Conley-Zehnder index of a periodic orbit to be well defined. 

Let $H: \Ss^1 
\times M \longrightarrow \RR$ be a time-dependent Hamiltonian and  
$X_H$ the Hamiltonian vector field defined by  $\iota_{X_H}\omega = 
dH(t,\cdot)$. Let $J_t \in \tx{End}(TM)$, $t\in \Ss^1$  
be a loop of almost complex structures which are compatible with the
 symplectic form i.e. 
$g_{J_t}(\cdot,\cdot) = \omega(\cdot, \, J_t \cdot )$, $t\in \Ss^1$  
are symmetric positive definite bilinear forms. We then have 
$J_tX_H(t,\cdot)=\nabla ^{g_{J_t}}H(t,\cdot)$. 
 
From a formal point of view, 
the construction of Floer homology closely follows that of Morse
homology with the following analogies:
\begin{eqnarray}   
\tx{Manifold } L &  \longleftrightarrow & 
\tx{Space of  contractible loops in } M, \tx{ denoted } \Lambda 
\nonumber \\  
& & \nonumber \\  
\tx{Metric } g &  \longleftrightarrow &  
L^2 \tx{ metric on } \Lambda~: \, \langle \eta, \, \xi 
\rangle = \int_{\Ss^1} g_{J_t} \big( \eta(t), \, \xi(t)  
\big) dt, \nonumber \\  
& & \hspace{5cm}  
\eta,\, \xi \in T_\gamma \Lambda= C^\infty 
\big(\Ss^1, \, \gamma ^* TM\big) \nonumber 
\end{eqnarray} 

\begin{eqnarray} 
\tx{Function } f &  \longleftrightarrow & \tx{Action functional }  
A_H~: \Lambda
 \longrightarrow \RR, \nonumber \\  
& &  \gamma \longmapsto - \int_{D^2} \overline 
\gamma ^* \omega - \int_{\Ss^1 }H(t,\gamma(t)) dt \nonumber \\
& & \nonumber \\
df(x)  &  \longleftrightarrow & dA_H (\gamma) \cdot \eta = 
- \int_{\Ss^1} \omega \big( \eta(t), \, \dot\gamma(t) - X_H(\gamma(t)) 
\big) dt \nonumber \\
& & \nonumber \\     
\tx{Critical point of } f  &  \longleftrightarrow &  
1-\tx{periodic orbit: } \dot \gamma(t) = X_H(t, \, \gamma(t)) 
\nonumber \\ 
&  & \nonumber \\       
\nabla ^g  f & \longleftrightarrow & \nabla ^J A_H (\gamma) = 
J_t\dot\gamma - J_tX_H \nonumber \\
\tx{Negative gradient trajectory }  & 
\longleftrightarrow  &  \tx{Map } u: \RR \longrightarrow 
\Lambda\ \tx{ or } \ u: \RR \times \Ss^1 
\longrightarrow M    \nonumber \\ 
 & & \tx{satisfying }  \nonumber \\  
& & \DP{u}{s} + J_t(u(s,t))  
\DP{u}{t}- \nabla ^{g_{J_t}} H \big(t,  u(s,t) \big) = 0 
\label{equation Floer} \\ 
& & \tx{and having bounded energy }  \nonumber \\  
& & E(u) = \frac12 \int_{\Ss^1}\int_{\RR}\Big( \big| \DP{u}{s}\big| ^2 
+ \big| \DP{u}{t} - X_H(t,u) \big|^2  \Big) \, ds \, dt < \infty  
\label{energie trajectoire Floer} 
\end{eqnarray} 
 
\medskip 

Let us make some comments before going on with the dictionary. For a
contractible loop  $\gamma: \Ss^1 
\longrightarrow M$ we have denoted  $\overline\gamma : D^2 
\longrightarrow M $ a map verifying   $\overline \gamma | _{ 
  \Ss^1} = \gamma $. The symplectic asphericity implies that
the action functional is well defined i.e. the quantity $\int_{D^2}
\overline \gamma ^*\omega$ does not depend on the choice of the
extension  $\overline \gamma$.  
 
From a formal point of view, the analogues of the negative gradient
trajectories of the finite dimensional case are maps 
 $u: \RR \times \Ss^1 
\longrightarrow M$ only verifying equation  (\ref{equation 
  Floer}). For finite-dimensional closed manifolds this
automatically implies the convergence at $\pm \infty$ to critical
points of $f$. On the contrary, when the dimension is infinite one has
to impose a supplementary condition in order to ensure the convergence
 at $\pm \infty$ to $1$-periodic orbits. 
This condition is precisely given by the finiteness of the energy
 (\ref{energie trajectoire
  Floer}). The latter is in fact  equivalent to the convergence at $\pm
\infty$ of a cylinder verifying (\ref{equation Floer}). 
 
The minima of the energy
are precisely the solutions of Floer's equation (\ref{equation Floer}).
For any two  orbits $x$, $y$ of period $1$ let us denote by
$\mathcal{U} (y, \, x)$ the space of maps  $u: \RR \times \Ss^1 
\longrightarrow M $  
that verify  $\lim_{s \rightarrow -\infty } u(s,\, t) 
=y(t)$, $\lim _{s \rightarrow +\infty} u(s,\, t) =x(t)$ with uniform
 convergence in  $t$. Any element in  $u \in \mathcal{U} (y, \, x)$ 
satisfies   
 
$$E(u) =  
\frac12 \int_{\Ss^1}\int_{\RR}\Big| \DP{u}{s} 
+ J_t(u) \DP{u}{t} - \nabla H(t,u) \Big|^2   ds \, dt + A_H(y) - 
A_H(x) \ .$$ 
 
\noindent  
This is a consequence of the identities
$- \frac{d}{ds} \, A_H(u(s,\cdot)) = \int_{\Ss^1} \omega(u_s, \, u_t - 
X_H) \, dt$ and 
$|u_s +Ju_t-JX_H|^2= |u_s|^2 +|u_t -X_H|^2 - 2 \int_{\Ss^1} 
\omega(u_s, \, u_t - X_H) \, dt$. The  minimum $A_H(y) - 
A_H(x)$ 
of the energy on  $\mathcal{U}(y, \, x)$  
is attained precisely by the solutions of equation  (\ref{equation
  Floer}) and we denote  $\mathcal{M}(y, \, x)$ 
the set of these trajectories. An
element  of
$\mathcal{M}(y, \, x)$ will be  called  a {\it Floer trajectory}.

A similar phenomenon occurs in the finite dimensional
case: the negative gradient trajectories are the minima of the energy 
 $E_{f,\, g}(\gamma) = \frac12 \int_\RR  
|\dot\gamma |^2 + |\nabla ^g f (\gamma) | ^2 \ ds$, defined on the
space of maps $\gamma: \RR \longrightarrow L$ that converge at
$\pm\infty$ to some given critical points of $f$. The transversality
condition that was mentioned in  \S\ref{subsection 1} can be rephrased
in terms of a surjectivity property for the linearization of the
equation of gradient trajectories \cite{Salamon Morse Floer} \, :  
$$D_\gamma \xi = \nabla _s \xi + \nabla_\xi \nabla f(\gamma), \qquad 
\xi \in \gamma ^* TL \ .$$ 
A similar transversality condition has to be verified in the infinite
dimensional case in order for the spaces  $\mathcal{M}(y, \, x)$ 
to inherit the structure of a finite dimensional manifold. The
linearization   
$$D_u \xi = \nabla _s \xi + J(u)\nabla_t\xi + \nabla_\xi J(u) u_t - 
\nabla _\xi \nabla H(t,u), \qquad \xi \in u^*TM $$ 
of the equation (\ref{equation Floer}) has to be a Fredholm map
whose differential is surjective at solutions of (\ref{equation
  Floer}). 
The Fredholm character requires the use of suitable functional spaces
for the elliptic analysis and  holds if the  $1$-periodic orbits of $H$
are nondegenerate. The  surjectivity is satisfied for a generic
choice of $H$ and $J$ \cite{FHS}. 
Our dictionary goes now as follows: 
\begin{eqnarray} 
\tx{Nondegenerate critical points} & 
  \longleftrightarrow &  1-\tx{periodic nondegenerate orbits } 
\nonumber \\  
& & \nonumber \\  
\tx{Morse index } \tx{ind}_{\tx{Morse}}(x)  
& \longleftrightarrow & \tx{Conley-Zehnder index } 
i_{CZ}(\gamma) \tx{ of a} \nonumber \\  
\tx{ of a critical point } x \qquad & & \tx{ periodic orbit 
  } \gamma, \text{\it taken with negative sign} \nonumber \\  
& & \nonumber \\ 
\dim \mathcal M (y, \, x) = \tx{ind}_{\tx{Morse}}(y) - 
  \tx{ind}_{\tx{Morse}}(x) & 
  \longleftrightarrow &  
 \begin{array}{rcl} & & \\
\dim \mathcal M (y, \, x) & = & i_{CZ}(x) -  i_{CZ}(y) \\
 & = & -  i_{CZ}(y) - \big( - i_{CZ}(x) \big)
\end{array} \nonumber 
\end{eqnarray} 
 
  Let us remark that, unlike in the finite dimensional case, the Hessian
  of $A_H$ at a critical point admits an infinite number of negative
  as well as positive eigenvalues. Its ``Morse'' index is therefore not
  well defined. This can be easily seen on the following
  example: consider $H\equiv 0$ on  $\CC$ and the action of a loop 
  $\gamma=\sum_{k \in \ZZ} z_k e 
  ^{ikt}$ is $\int xdy= \pi \sum _{k \in \ZZ} k |z_k|^2$. The positive
  and negative eigenspaces of this quadratic form are obviously
  infinite dimensional. In fact, this shows again that the classical
  Morse theory is not adapted to the present context: it would come to
  glueing infinite dimensional discs along infinite dimensional
  spheres, which are contractible. The homotopical invariants of the
  corresponding spaces would vanish at any step of the construction. 
 
  The Conley-Zehnder index is an integer that is associated to a path
  of symplectic matrices having the identity as  origin and whose end
  does not contain the eigenvalue $1$ in its spectrum. To any periodic
  orbit 
  $\gamma$ one can associate a  Conley-Zehnder index by trivializing
  $TM$ over a filling disc $\overline \gamma$ and by considering the
  path of symplectic matrices induced by the linearization of the
  Hamiltonian flow along $\gamma$. The assumption $\langle c_1, \,
  \pi_2(M) \rangle =0$ ensures that two such trivializations are
  homotopically equivalent along $\gamma$ and the integer 
  $i_{CZ}(\gamma)$ will thus be independent of the trivialization. 
  We explain at the end of this section why the suitable analogue of
  the Morse index is the Conley-Zehnder index considered with negative
  sign, rather than simply the  Conley-Zehnder index. 

  If the transversality assumptions are verified, the implicit
  function theorem ensures that the dimension of 
  $\mathcal M (y, \, x)$ at $u$ is equal to the Fredholm index of 
  $D_u$. The identification of the latter with the difference between
  the Conley-Zehnder indices of the ends is a consequence of a
  characterization in terms of the spectral flow of a certain family
  of first order differential operators that is associated to equation 
  (\ref{equation Floer}) \cite{RS1, RS2, Sa}.  
 
  The only ingredient still lacking in order to formally define a
  (co)homological differential complex in analogy with the TSW complex
  is a recipe to associate a sign to a Floer trajectory running
  between periodic orbits whose Conley-Zehnder indices have 
  a difference
  equal to $1$. We shall not pursue this matter here and will just
  claim  that {\it there is} one such recipe
 \cite{Floer Hofer orientations, Salamon Morse Floer}. One can
 alternatively work with $\ZZ/2\ZZ$ coefficients in order to avoid
 all sign problems.
 
\bigskip 

\noindent {\sf \large Definition.} {\it The homological Floer complex 
  $FC_*(M; \, H, \, J)$   
 and the cohomological one\break $FC^*(M; \, H, \, J)$ are 
 defined respectively by the analogues of equations
(\thesection.\thefirstTSW, \thesection.\thesecondTSW,
\thesection.\thethirdTSW) and     
(\thesection.\thefirstTSW$^\prime$, \thesection.\thesecondTSW$^\prime$,
\thesection.\thethirdTSW$^\prime$) through the above dictionary. 
The grading is given by the opposite of the
 Conley-Zehnder index. }  

\bigskip 
 
The definition depends on the transversality 
results that we have mentioned above. In the case of a Hamiltonian
having nondegenerate periodic orbits, the almost complex structures
for which these  hold form a set 
$\mathcal{J}_{reg}(H)$ which is of the second category in the sense of
Baire in the space of $\omega$-compatible almost complex structures. 
Conversely, for any fixed family  $J=(J_t)$, 
the transversality results are valid for a second Baire category set
of Hamiltonians  $\mathcal H _{reg}(J) \subset 
C^{\infty}(\Ss^1 \times M, \, \RR)$. In finite dimension, this amounts
to prescribe the metric and to choose a generic Morse function: from
the point of view of genericity, the metric and the function
play symmetric roles. 

In order to fix ideas, all the remarks that are to follow will concern
cohomology groups. The distinction homology - cohomology will gain
importance only for manifolds with boundary. The definition will in
that case 
contain as a supplementary ingredient an algebraic limit process
which, according to the formalism being homological or cohomological,
will be direct or inverse.

The fundamental property of Floer cohomology
$$FH^*(M; \, H, \, J) \, = \, H^*(FC^*(M; \, H, \, 
J))$$   
is its independence with respect to the Hamiltonian and with respect
to the almost complex structure. For any two pairs
$(H^0, \, 
J^0)$, $(H^1, \, J^1)$ which satisfy the above regularity conditions
there is a homotopy of regular pairs $(H^t, \, J^t)$, $t \in [0, \,
1]$ that links them together. Any such homotopy induces an {\it
  isomorphism} $FH^*(M; \, H^0, \, J^0) 
\stackrel\sim\longrightarrow  FH^*(M; \, H^1, \, J^1)$ which, moreover,  
{\it does not depend on the chosen regular homotopy.} 
The consequence is  that one can identify Floer (co)homology
with a classical topological invariant, namely singular (co)homology. 
Let us recall the relevant argument. 
Consider a Hamiltonian function that is time
independent and Morse. Any critical point $x$ of $H$ is a (constant)
$1$-periodic orbit and, in view of the convention
$X_H = -J \nabla H$, we infer 
$i_{CZ}(x) = n - \tx{ind}_{\tx{Morse}}(x,
\, -H)$, $n= \frac 1 2 \dim M$ or, written in a different way,
 $\tx{ind}_{\tx{Morse}}(x, \,
-H) = n + \big( - i_{CZ}(x)   \big)$. When  $H$ is  small enough
in the  $C^2$ norm one can show that  there are no $1$-periodic orbits 
other than the critical points of $H$ and, moreover, the 
Floer trajectories 
(solving  $u_s + Ju_t = \nabla H$) 
that run between points whose index difference is equal to $1$ are in
fact time independent \cite{F2}. This means that the Floer complex
coincides with the TSW complex corresponding to  $\nabla H$. If 
the grading on  $FC^*$ is given by {\it minus} the  Conley-Zehnder
index, one gets the isomorphism 
$$FH^*(M, \, \omega) \simeq H^{n+*}(M; \, \ZZ), \quad 
n= \frac 1 2 \dim M \ .$$ 
We see in particular the interest of  grading by
$-i_{CZ}$ rather than  $i_{CZ}$. In the latter case we would have
obtained the isomorphism  $FH^* \simeq H^{n+*}$ through the extra use
of  Poincar\'e duality $H_{n-*}(M) \simeq H^{n+*}(M)$.

\medskip

\subsection{Floer homology for manifolds with contact type boundary or
Symplectic homology} 
\label{hom Floer type contact}  
The second big setting in which Hamiltonian Floer homology groups can
be defined is that of {\it compact symplectic manifolds with contact
  type boundary}. The main references on this topic are the papers by 
K. Cieliebak, 
A. Floer, H. Hofer, K. Wysocki \cite{FH, CFH, CFHW} and C. Viterbo 
\cite{functors1}. I have used D. Hermann's thesis \cite{these Hermann} 
with great profit due to the very clear exposition style. The paper of 
P. Biran, L. Polterovich and D. Salamon \cite{BPS} contains results
concerning the existence of periodic orbits representing nontrivial free
homotopy classes. The  
Weinstein conjecture and related problems are discussed in the book by
 H. Hofer and  E. Zehnder \cite{HZ}. 
 
The initial motivation for the construction of Floer homology groups was
the  existence problem for closed orbits of Hamiltonian
systems. Roughly speaking, this problem has  
two distinct aspects: existence of closed orbits with a {\it given
  period} and existence of closed orbits on a {\it given energy
  level}. Historically, these two directions correspond to two
conjectures of  V.I. Arnold \cite{Ar} (1965) and 
A. Weinstein \cite{We}  
(1979). The first claims that a lower bound for the number of closed orbits
with fixed period on a closed symplectic manifold $M$ is provided by the
rank of the (rational) cohomology $H^*(M; \, \QQ)$. Under the
assumptions of the previous section, this follows from the very
construction of the Floer homology groups. Since 1996 there are proofs
that work in full generality, with no extra assumptions on the
underlying manifold. The second conjecture claims the existence of at
least one closed orbit on a regular compact energy level $\Sigma$
which is of contact type. Floer homology for manifolds with boundary
is a tool that is particularly adapted to the study of this question
(\S\ref{Weinstein}), but we shall  present two other fascinating 
 applications  in \S\ref{poly ell} and \S\ref{stab}. 

\begin{defi} A  (compact) hypersurface $\Sigma$ of a symplectic
  manifold 
  $(M, \, \omega)$ is said to be  {\rm of contact type} if there is a
  vector field $X$ defined in a neighbourhood of $\Sigma$, transverse
  to $\Sigma$ and verifying  $L_X\omega = \omega$. The vector field  
  $X$ is called  the {\rm  Liouville field.} The $1$-form $\lambda = 
  \iota _X\omega $ is called the  {\rm Liouville form.} If the
  Liouville field is globally defined on the whole of $M$ we say that  
  $\Sigma$ is of  {\rm restricted contact type}.  
 
  The boundary of a compact symplectic manifold $M$ is said to be of
  {\rm (restricted) contact type} if the above conditions are
  satisfied and the Liouville field is outward pointing. 
\end{defi} 
 
The contact type condition is a symplectic analogue for convexity in the
linear symplectic space $\RR^{2n}$: any compact convex  hypersurface is of
(restricted) contact type, as the radial vector field
$X(x)=\frac12 x$, $x \in \RR^{2n}$  
satisfies the above conditions (assuming $0$ is in the bounded component
of   $\RR^{2n} \setminus \Sigma$).  
The conjecture has been formulated precisely in view of preliminary
existence results  on convex or star-shaped energy
levels. A first proof for a contact type $\Sigma  
\subset \RR^{2n}$ was given by  C. Viterbo \cite{IHP}
and lots of other ambient spaces have been subsequently exlored.

The contact type condition is related to holomorphic pseudo-convexity,
as remarked by 
Y. Eliashberg, M. Gromov \cite{EG} and D. McDuff \cite{dusa}. This is
precisely the reason why we impose that the Liouville field be outward
pointing. One should 
remark that this is automatically true if the boundary is of 
restricted
contact type as the Liouville field (exponentially) expands
volumes. We give more details on pseudo-convexity
in the sequel.

Floer homology groups of a manifold $M$ with contact type boundary
will be
defined with the help of Hamiltonians that admit the boundary
$\partial M$ as a regular level. The invariants that we thus obtain
will take into account not only the $1$-periodic orbits in the
interior of the manifold, but also the closed orbits having arbitrary
period on the boundary. As a consequence, they are well adapted to
the study of Weinstein's problem. Let us stress from the very
beginning that, unlike Floer homology of closed manifolds which  
is finally proved to be equal to the singular homology,
Floer homology of manifolds with boundary has no similar topological
correspondent. It is all the more important  to
exhibit in the latter case 
qualitative properties that are
determined by  additional geometric properties of the manifold. 

\medskip Here is how one retrieves closed orbits of arbitrary period on
a contact type level with the help of $1$-periodic orbits in a
neighbourhood of $\Sigma$. Note  that the restriction of 
$\omega$ to  $T\Sigma$ has a one dimensional kernel on which 
$\lambda$ does not vanish. If $H$ is an autonomous Hamiltonian
admitting  
$\Sigma$ as a regular level then $X_H \in \ker \, \omega | 
_{T\Sigma}$ and $\lambda(X_H) \neq 0$.  
 
\begin{defi} \label{defi Reeb} 
The {\rm  Reeb vector field (or characteristic field)}
  $X_{\tx{Reeb}}$   
of $\Sigma$ is defined by the following two properties:
$X_{\tx{Reeb}} \in \ker \, \omega |_{T\Sigma}$ and 
$\lambda(X_{\tx{Reeb}})=1$. An orbit of  $X_{\tx{Reeb}}$ is called a  {\rm 
  characteristic of $\Sigma$.}  
\end{defi}  
 
\noindent One should note that the area  $\int_\gamma \lambda$  
of a closed characteristic is equal to its period. Moreover, the
orbits of $X_H$ that are located on $\Sigma$ are in one-to-one correspondence
with the latter's characteristics. It is important to understand that
the Hamiltonian dynamics on a regular level does not depend on the
Hamiltonian but on the level itself: it is more of a geometric nature
 rather than purely dynamic or analytic.

Let us denote by  $\varphi ^t_X$ the flow of the Liouville field. A
whole neighbourhood  $\mc{V}$ of $\Sigma$ is foliated by the
hypersurfaces  
$\varphi ^t_X(\Sigma) _{-\delta < \, t < \, \delta}$ with $\delta > 0$ 
small enough. In view of ${\varphi ^{t \ *}_X }  \omega = e ^t 
\omega$, the characteristics on these hypersurfaces are in one to one
correspondence with those of $\Sigma$. It is now comfortable to make a
coordinate change via the symplectic diffeomorphism
\begin{equation} \label{diffeo Liouville}  
\Psi~: \Sigma \times [1- \delta, \, 1+\delta] \stackrel \sim 
  \longrightarrow  
  \mc{V}, \qquad \delta > 0 \tx{ small} \ ,  
\end{equation}  
  $$\Psi (p,\, S) = \varphi ^{\ln(S)}_X(p) \ ,$$ 
  verifying
  $$\Psi ^* \lambda = S\cdot \lambda |  $$  
  where $\lambda |$ is the restriction of  $\lambda= \iota_X\omega$  
to $T\Sigma$. The autonomous Hamiltonian  
$$\begin{array}{l}  
H~: \Sigma \times [1- \delta, \, 1+\delta] \longrightarrow \RR \\  
{} \\ 
H(p,\, S) = h(S), \qquad h~: [1- \delta, \, 1+\delta] \longrightarrow \RR 
\end{array} 
$$ 
satisfies $X_H(p,\, S) = -h'(S)X_{\tx{Reeb}}$. Its $1$-periodic orbits
that are located on the level  $S$ correspond to characteristics of
period $h'(S)$ located on  $\Sigma$ (with the opposite orientation). 
By studying the  $1$-periodic orbits of such Hamiltonians one will have
in fact studied characteristics on $\Sigma$: the more important the
variation of  $h$ in the small interval  $[1-\delta, \, 
1+\delta]$ is, the more characteristics one ``sees''. 

\medskip

Any reasonable Floer homology invariant for a manifold with boundary
should take into account the topology of the manifold and {\it all}
closed characteristics on the boundary. We thus retrieve the common
underlying principle of the constructions in 
\cite{FH, CFH, CFHW, functors1}: the cohomology  
\begin{equation} \label{homologie Floer et limite directe}  
\displaystyle FH^*(M)= 
\lim_{(H, \, J)} FH^* (M; \, H, \, J) 
\end{equation}   
will be defined as a limit following an {\it admissible family} of
Hamiltonians, steeper and steeper near the boundary.  As we have
already warned the reader, this limit is of direct or inverse type
respectively, according to the choice of a homological or
cohomological formalism. One supplementary refinement will consist in 
using a
truncation by the values of the action, a crucial ingredient for
the applications in \S\ref{poly ell} and \S\ref{stab}. Before reaching
them, let us  describe  the various points of view already present
in the literature. 
 
\bigskip 
 
\subsubsection{\sf \large Symplectic homology of a bounded open set   
$U \subset  (\CC^n, \, \omega_0)$, cf. \cite{FH}}   
\label{homologie FH} 
This definition was introduced by A. Floer and 
  H. Hofer \cite{FH}. In order to give a presentation in tune with the
  subsequent constructions  
we shall present below the cohomological
  setup, while the original point of view is homological. The two main
  features of this theory are the following:
 \renewcommand{\theenumi}{\alph{enumi}} 
  \begin{enumerate}  
    \item it is ``extrinsic'' in the sense that Hamiltonians are
    defined on the whole ambient space $\CC^n$; 
    \item the definition is valid for arbitrary open sets, without any
    regularity or contact type hypothesis on the boundary. 
  \end{enumerate}  

The class  $\mathcal{H}(U)$ of admissible Hamiltonians  $H : 
  \Ss^1 \times \CC ^n \longrightarrow \RR$ is defined by the following
  properties:
\renewcommand{\theenumi}{\arabic{enumi}} 
\begin{enumerate}  
   \item \label{one} $H_{\vert_{ \Ss^1 \times \bar{U}}} <0$;  
   \item \label{two} there is a positive definite matrix  $A$  such that   
        $$ \frac{\vert H'(t,u) - Au \vert  }{\vert u\vert }   
        {\longrightarrow} 0, \qquad \vert u \vert \rightarrow  \infty $$  
   uniformly in  $t \in \Ss^1$;  
   \item \label{three} the differential system    
        $$ -i\dot{x}= Ax , \qquad x(0) = x(1)$$  
    admits only the trivial solution  $x \equiv 0$;  
   \item \label{four} there is a constant  $c > 0$ such that    
      $$\parallel H''(t,u)\parallel  \leq  c , \qquad \forall \ t \in 
      \Ss^1, \ u \in \CC^n \ ,$$  
      $$ \Big| \DP{H'}{t} (t,u) \Big| \leq c(1 + \vert u \vert ),  
      \qquad \forall \ t \in \Ss^1, \ u \in \CC^n \ .$$  
  \end{enumerate}  

Condition (\ref{one}) ``pins down'' the Hamiltonian along $U$ and
allows it to (steeply) increase only in the neighbourhood of the
boundary, according to the philosophy that we have detailed above. Condition
(\ref{two}) prescribes the asymptotic behaviour of $X_H$ and, combined
with (\ref{three}), ensures not only that all its $1$-periodic
orbits are contained in a compact set, but also that the Floer
trajectories that link them together do not escape to infinity. This
is a crucial point in proving Floer compactness for open manifolds and
we shall give more details on it in the context of 
the two  constructions to follow. 
Condition (\ref{four}) is of a  technical nature and plays some 
role in the proof of the $C^0$ estimates.  
One should note that condition  (\ref{two}) is to be interpreted as  a
kind of quadratic asymptoticity  for
$H$. 

We denote by $\mathcal{H}_{reg}(U)$ the class of admissible
Hamiltonians having nondegenerate $1$-periodic orbits: it is of 
second Baire category in  $\mathcal{H}(U)$. We denote by
  $\mathcal{J}$ the class of almost complex structures on $\CC^n$ that
  are compatible with  $\omega_0$ and that are equal to the standard
  complex structure $i$ outside a compact set. The transversality in
  Floer's equation is verified for a  dense set 
$\mathcal{HJ}_{reg}(U) \subset \mathcal{H}_{reg}(U) \times 
\mathcal{J}$ and the truncated Floer cohomology groups are defined for
a regular pair as follows: 
\begin{eqnarray*}
& & FC^k_{]a, \, +\infty[}(H, \, J) =  
          \bigoplus_{ {\tiny \begin{array}{c} 
                                - i_{CZ}(x) = k \\  
                                A_H(x) > a  
                             \end{array} } }  
           \ZZ \langle x \rangle ,  
           \qquad a\in \RR \cup \{ -\infty  \} \ , \\  
& & 
FH^*_{]a, \, +\infty[}(H, \, J) = H^* (FC^*_{]a,\, +\infty[} (H, \,
           J)) \ , 
\end{eqnarray*} 
\begin{eqnarray*}
& & FC^*_{]a,b]} (H, \, J) =  
FC^*_{]a, \, +\infty[}(H, \, J) \ / \ FC^*_{]b, \, +\infty[}(H, \, J),  
\qquad -\infty \le a < b < +\infty \ , \\ 
& & \\
& & FH_{]a,b]}^*(H, \, J)=H^*(FC^*_{]a,b]} (H, \, J)) \ .
\end{eqnarray*} 
The inverse limit 
 (\ref{homologie Floer et limite 
  directe})  is considered with respect to the following partial order
 relation on  $\mathcal{HJ}_{reg}(U)$, which induces an inverse
 system on the cohomology groups: 
$$(H, \, J) \prec (K, \, \widetilde{J}) \qquad \tx{iff} \qquad H(t,\, u) 
\le K(t,\, u) \ .$$ 
The cohomological inverse system is determined in the following
manner. For two ordered pairs
$(H, \, J) \prec (K, \, \widetilde{J})$ 
we consider a homotopy  $(H(s, \, t, \, u), \, J(s, \, t, \, 
u))$ such that:  
\vspace{-.3cm}  
\begin{itemize}  
\item There is an $s_0 >0$ with   
$(H(s, \, t, \, u), \, J(s, \, t, \, 
u)) \equiv \begin{array}{ll} {} & {} \\  
(H(t, \, u), \, J(t, \, u)), & s \le -s_0 
  \\ (K(t, \, u), \, \widetilde{J}(t, \, u)), & s \ge s_0 \quad ;  
\end{array}$  
\item The homotopy is an increasing function with respect to $s$ i.e.  
$\frac{\partial H}{\partial s} \, (s, \, t, 
\, u) \ge 0$;  
\item The homotopy 
satisfies certain additional assumptions concerning the
  behaviour at infinity and the regularity. The most important
  of them is the following compatibility with the quadratic
  asymptoticity: 

\qquad There is a smooth path $A(s) \in \tx{End}_{\, \RR}(\CC^n)$ 
of positive definite matrices such that $A(s)=A(-s_0)$ for 
$s\le -s_0$, $A(s) = A(s_0)$ for $s\ge s_0$ and 
$$\frac {|H'(s,\, t, \, u) - A(s)u|}{|u|} \longrightarrow 0, \qquad
|u| \rightarrow \infty \ .$$
Moreover, we require that if the Hamiltonian system $-i \dot{x}=
A(\hat{s})x$ has a nontrivial $1$-periodic solution for some
$\hat{s}\in \RR$, then $\frac d {ds}A(s)\Big|_{s=\hat{s}}$ is positive
definite. 
\end{itemize}  
The asymptotic conditions on the homotopy ensure that the solutions of
the {\it parameterized Floer equation}
\begin{equation} \label{Fl param 1} 
u_s + J(s, \, t, \, u) u_t - \nabla \, H(s, \, t , \, 
u) =0 \ , 
\end{equation}  
\begin{equation} \label{Fl param 2}  
u(s, \, t) \longrightarrow x^\pm, \ s \longrightarrow \pm \infty \ , 
\end{equation} 
which play now the role of the trajectories verifying
(\ref{equation Floer} - \ref{energie trajectoire Floer}), do also stay
in a compact set. Here  $x^-$ and $x^+$ are 
$1$-periodic orbits of  $H$ and  $K$ 
respectively. The dimension of the moduli 
space $\mc{M}(x^-, \, x^+)$ of solutions of   
(\ref{Fl param 1} - \ref{Fl param 2}) is  
$-i_{CZ}(x^-)- \big( - i_{CZ}(x^+) \big)$ 
but, to the difference of equation (\ref{equation 
  Floer} - \ref{energie trajectoire Floer}), the group 
$\RR$ is no longer acting by translation. This implies that,
generically,  the moduli space is no longer empty if its formal
dimension is zero and the morphism of complexes defined according to
(\thesection.\thethirdTSW$^\prime$) will  respect the
grading: 
\begin{equation} \label{syst proj} 
 \sigma~: FC ^*_{]a, \, \infty]}(K, \, \widetilde{J}) \longrightarrow 
FC^*_{]a, \, \infty]}(H, \, J) \ , 
\end{equation}  
$$\sigma \langle x^+ \rangle = \sum_{i_{CZ}(x^-)=i_{CZ}(x^+)} \# 
\mc{M}(x^-, \, x^+) \, \langle x^- \rangle \ .$$ 
The map  $\sigma$ is well defined because the action is decreasing
along solutions of the parameterized equation 
 (\ref{Fl param 1}). This is due to the fact that the homotopies that
 we allow are increasing:
\begin{eqnarray*} \frac{d}{ds} \, A_{H(s)}(u(s,\cdot)) & = &  
- \int_{\Ss^1} \omega(u_s, \, u_t - X_{H(s)}(u(s, \, t))) \, dt  -  
\int_{\Ss^1} \frac{\partial H}{\partial s}(s, \, t, \, u(s, \, t)) \, 
dt \\ 
& = & -\int_{\Ss^1} \omega ( u_s, \, J(s, \, t, \, u(s,t)) u_s ) \, dt 
- \int_{\Ss^1} \frac{\partial H}{\partial s}(s, \, t, \, u(s, \, t)) 
\, dt \\  
& = & -\parallel u_s \parallel  _{g_{_{J(s)}}} ^2 -  
\int_{\Ss^1} \frac{\partial H}{\partial s}(s, \, t, \, u(s, \, t)) 
\, dt  \ \le \ 0 \ .  
\end{eqnarray*}  
The morphism $\sigma$ commutes with the differential and descends to a
morphism in cohomology that is called ``monotonicity map''
\begin{equation} \label{et voila le morphisme de monotonie}  
\xymatrix{ 
FH^*_{]a,b]}(K, \, \widetilde{J})  
\ar[rr]^{\ \sigma_{(H, \, J)}^{(K, \, \widetilde{J})}} & &  
FH^*_{]a,b]}(H, \, J) \ . 
}
\end{equation} 
As in the compact case, two admissible homotopies induce the same
morphism in homology and we moreover have 
$$\sigma_{(H, \, J)}^{(K, \, \widetilde{J})}  
\ \circ \   
\sigma ^{(K', \, \widetilde{J}')}_{(K, \, \widetilde{J})}  
\ = \ 
\sigma ^{(K', \, \widetilde{J}')}_{(H, \, J)} \ , $$ 
for $(H, \, J) \prec (K, \, \widetilde{J}) \prec (K', \, 
\widetilde{J}')$. Let us then define   
\begin{eqnarray*}
& & \displaystyle FH^*_{]a,b]}(U)= 
\lim_{\stackrel{\longleftarrow}{(H, \, J)}}   
      FH_{]a,b]}^*(H, \, J) \ , \qquad -\infty \le a < b < +\infty \ ,
      \\ 
& & \displaystyle FH^*_{]a, \, +\infty[}(U)= 
\lim_{\stackrel{\longleftarrow}{(H, \, J)}}   
      FH_{]a,\, +\infty[}^*(H, \, J) 
\end{eqnarray*}
and get in particular
$$\displaystyle FH^*_{]a, \, +\infty[}(U)= \lim_{\stackrel \longleftarrow b}  
FH_{]a,b]}^*(U), \qquad b\rightarrow +\infty \ .$$    
The ``truncation maps'' $FH^*_{]a, \, b]}(U) \longrightarrow 
     FH^*_{]a, \, b']}(U) $, $b \ge b'$ are induced by the inclusions 
     $FC^*_{]a, \, \infty[} (H, \, J) \hookrightarrow 
     FC^*_{]a', \, \infty[}(H, \, J)$, $a \ge a'$ 
which give rise to the morphisms 
$$FH^*_{]a, \, b]}(H, \, J) \longrightarrow FH^*_{]a', \, b']}(H, 
     \, J), \qquad a \ge a', \ b \ge b' \, .$$  
The latter are compatible with the monotonicity maps. 
For a fixed value of  $a$ the groups  
     $FH^*_{]a, \, b]}(H, \, J)$ form an inverse bi-directed
     system. In view of the fact that two inverse (or direct) limits
     in a bi-directed system commute one can therefore write 
     $$FH^*_{]a, \, +\infty[}(U) 
= \lim_{\stackrel \longleftarrow b} \lim_{\stackrel 
     \longleftarrow {(H, \, J)}} FH^*_{]a, \, b]}(H, \, J) =  
\lim_{\stackrel 
     \longleftarrow {(H, \, J)}} \lim_{\stackrel \longleftarrow b}  
FH^*_{]a, \, b]}(H, \, J), \qquad  b \rightarrow \infty \ . $$ 
We shall see below  (cf. \S\ref{relations between homologies}) 
that, for an open set  $U$ with restricted contact type boundary, the
groups  $FH^*_{]a, \, +\infty[}(U)$ do not depend on  $a$ 
if the latter is strictly negative. This is due to the
existence of a cofinal family whose  $1$-periodic orbits all have
an action that is either positive, either negative and arbitrarily
close to zero 
\cite{these Hermann, functors1}. We define in the general case  
\begin{equation}
  FH^*(U) = FH^*_{]-\infty, \, +\infty[}(U) \ .
\end{equation} 
We shall then have 
$FH^*(U) = FH^*_{]a, \, +\infty[}(U)$, $a < 0$ if $U$ has restricted
contact type boundary. 

\medskip

As we have already pointed out, 
the present version of Floer homology asks for an important new
ingredient compared to the case of closed manifolds, 
namely the existence of a priori $C^0$ bounds 
which ensure that, for fixed limiting orbits, all Floer 
trajectories for the parameterized equation stay in a 
compact set. If this were not the case then some pathological 
noncompactness might appear in the moduli space of Floer trajectories: 
 the number of trajectories between two orbits 
with index difference equal to one might be infinite, or the square of
the Floer differential might no longer be zero. We give more details
about the proofs of these bounds for the two constructions to follow
and in  \S\ref{C 0 bounds}.

\medskip  
 
Let us mention that an analogous construction  can be used to define
Floer {\it homology} groups. The differential will be defined
according to  
 (\thesection.\thethirdTSW) (and this is the original setup of 
 \cite{FH}). The sub-complexes that will be  preserved by the
 differential are 
$$FC_*^{]-\infty, \, a]}(H, \, J) = \bigoplus_{A_H(x)<a}\ZZ \langle x 
\rangle \ , \qquad a \le \infty    
$$ 
and the monotonicity and truncation morphisms will define a {\it
  direct} double system: 
$$FH_*^{]a, \, +\infty[}(U) 
= \lim_{\stackrel \longrightarrow b} \lim_{\stackrel 
     \longrightarrow {(H, \, J)}} FH_*^{]a, \, b]}(H, \, J) =  
\lim_{\stackrel 
     \longrightarrow {(H, \, J)}} \lim_{\stackrel \longrightarrow b}  
FH^*_{]a, \, b]}(H, \, J), \qquad  b \rightarrow \infty \ . $$ 
 
This distinction bears a specific importance: the direct limit is an
exact functor, while the inverse limit is generally only left exact
and becomes exact if the terms of the directed system are all finite
dimensional vector spaces (\cite{ES} Ch. VIII).  
This has nontrivial consequences in practice: one can prove for
example a K\"unneth formula that is valid in Floer {\it homology} with
arbitrary coefficients, while the analogous K\"unneth formula holds in
Floer {\it cohomology} only with coefficients in a field (see
\S\ref{Weinstein}).

\medskip 

\subsubsection
{\sf \large Symplectic homology of relatively compact open sets in
  manifolds with contact type boundary, cf. \cite{CFH} and \cite{CFHW}.}  
 \label{CiFH}
    We describe below a construction introduced by  
    K. Cieliebak, A. Floer and H. Hofer in \cite{CFH}, as well as one
    of its variants appearing in \cite{CFHW}. 
 
         The admissible Hamiltonians for a relatively compact open set 
    $U \subset M \setminus \partial M$ are defined in \cite{CFH} 
   by the following
    properties: 
    \begin{itemize}  
     \item $H_{\vert_{\Ss^1 \times \bar{U}}} <0$;   
     \item $H \equiv m(H) = \max H >0$ in a neighbourhood of $\partial M$;  
     \item all periodic orbits satisfying  $\int_0^1 H(t,x(t)) dt <  
       m(H)$ are nondegenerate. 
    \end{itemize} 

    A neighbourhood of the boundary is trivialized by the Liouville
    flow $\varphi \circ \ln$ 
    as $\big( \partial M \times [1-\delta, \, 1], \, d(S\lambda|)
    \big)$, $\delta >0$. The
    admissible almost complex structures are defined to be those
    that can be written near the boundary as 

\begin{equation} \label{std a c str} 
\left\{ \begin{array}{l} J_{(p,S)} | _\xi = J_0 \, ,
\\ J_{(p,S)}(\DP{}{S}) = \frac 1 {CS} 
    X_{\tx{Reeb}} (p) \, , \quad C >0 \, ,
\\ J_{(p,S)}(X_{\tx{Reeb}}(p)) = -CS \DP{}{S}  \, .
\end{array}  \right. \, ,
\end{equation}  
with $J_0$ an almost complex structure compatible with the restriction
of  $\omega$ to the contact distribution $\xi$ on $\partial M$. 
    These are precisely the almost complex structures that are
    invariant by homotheties $(p, \, S) \longmapsto (p, \, aS)$,
    $a>0$. 
The homotopies are again chosen to be increasing in the $s$-variable
and such that, for fixed $s$, the Hamiltonian $H(s, \cdot, \cdot)$ is
of the type above. 

The definitions of the Floer complex and of the (co)homology
    groups are perfectly similar to those that we have described for
    open sets in $\CC^n$ and the existence of the  $C^0$ bounds for
    the parameterized trajectories is again the
    crucial point. In the present situation it is  replaced by
    the requirement that the parameterized Floer trajectories stay at
    bounded distance from the boundary. To the difference of the
    preceeding construction and thanks to the very special form of the
    Hamiltonians and of the homotopies, a geometric argument based on
    holomorphic convexity allows one to conclude easily.
 
\begin{defi}  
Let  $J$ an almost complex structure that is compatible  with the
symplectic form  $\omega$. A hypersurface $\Sigma \subset M$ is said
to be  
$J$-convex if it can be locally written as the regular level set of a
plurisubharmonic function i.e. a function $\varphi : M \longrightarrow
\RR $ which satisfies 
$dd^c\varphi \le 0$, where  $d^c = J^*d$. 
\end{defi}  
 
\noindent {\sf Example:} if   
$M$ is a manifold with contact type boundary, all hypersurfaces
  $\partial
M \times  
\{ S_0 \} \subset \partial M \times [1-\epsilon, \, 1]$ are 
 $J$-convex with respect to any admissible almost complex structure.  
Indeed, $\varphi(p,\, S) = S$ is plurisubharmonic as
    $dd^cS = d(-CS\lambda|)=-C\omega \le 0$. 
 
\medskip 

The interest of plurisubharmonic functions and pseudo-convex
hypersurfaces lies in the following lemma. 

\begin{lem} {\sf (H. Hopf, see also  \cite{dusa})} \label{max princ} 
Let  $\Sigma \subset M$ be a  
    $J$-convex hypersurface and  $\varphi$ a (local) function of
    definition. No  $J$-holomorphic curve  $u:(D^2(0, \, 1), \, i) 
    \longrightarrow M$ can have an interior strict 
    tangency point with  $\Sigma$  
    i.e.  $\varphi \circ u $ cannot have a strict local maximum.   
\end{lem}  
 
\demo Let  $z=s+it$ be the complex coordinate on  $D^2(0, \, 1)$ and 
 $J_0=i$ be the standard complex  structure 
on the disk.  The crucial (easy) computation is 
$dd^c_{J_0}(\varphi \circ u) = - \Delta(\varphi \circ u ) ds 
\wedge dt$. The  $J$-holomorphicity of  $u$ implies   
$dd^c_{J_0}(\varphi \circ u) = dJ_0^*u ^* d\varphi =  d u ^* J^* d 
\varphi = u ^* dd^c_J \varphi $. But  
$\varphi$ is plurisubharmonic and thus 
 $\Delta(\varphi \circ u) \ge 0$ i.e. $\varphi \circ u$ is
 plurisubharmonic in the classical sense. Hence $\varphi \circ u$ 
 satisfies the mean
 value inequality and there can be no strict maximum that is achieved
 in the interior of the disk. 

\hfill{$\square$} 

\medskip  

This shows that all Floer trajectories corresponding to a {\it fixed} 
Hamiltonian stay at a bounded distance from the boundary. The argument
requires some further refinement in order to deal with parameterized
trajectories but its basic feature remains the use of the maximum
principle \cite{CFH}.

    \medskip  
 
    In order to get an intuition of what the corresponding invariants
    really compute, let us consider the case where 
    the open set $U$ has a contact type boundary as well.
    One can then construct a cofinal family of Hamiltonians with a special
    geometric meaning, whose form is sketched 
    in Figure \ref{les trois types de familles admissibles}
    (2). A  neighbourhood of $\partial U$
    is trivialized by the corresponding Liouville flow as $\partial U
    \times [1-\epsilon, \, 1+\epsilon]$  
    and the Hamiltonians  are of the form $H=h(S')$, $S'\in
    [1-\epsilon, \, 1+\epsilon]$, with $h$ steeper and steeper (see
    \S\ref{Homologie symplectique Viterbo} for a discussion of
    nondegeneracy).  
    Characteristics on
    $\partial U$ are therefore 
    seen twice, once on the convex part and once on the concave part
    of the  function $h$.

\begin{figure}[h]  
      \begin{center}   
      \includegraphics{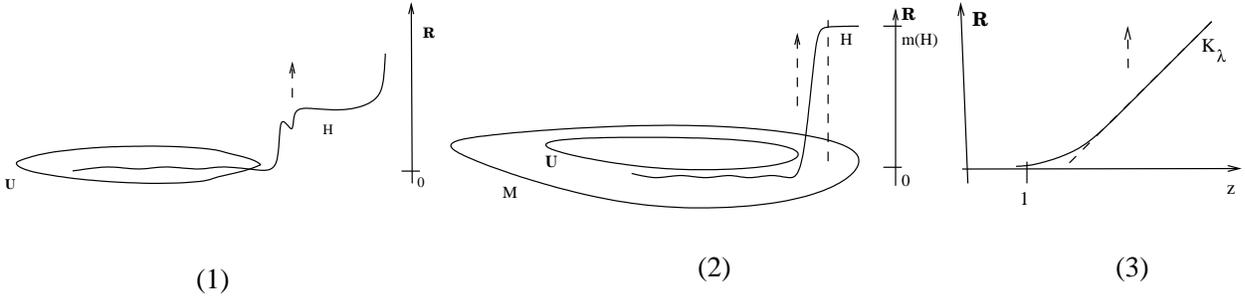}  
      \end{center}   
      \caption{{Admissible families: (1) - \cite{FH}, \ (2) - 
          \cite{CFH}, \ (3) - \cite{functors1}}\label{les trois types 
          de familles admissibles}}   
\end{figure}   

Here is now the slight modification of the above definition which is
  used in 
  \cite{CFHW} for  proving the 
  stability of  the action spectrum for a compact symplectic
  manifold with contact type boundary. We describe it
  in view of its use in \S\ref{stab}, as well as in order to further
  discuss the relevance of the direction in which morphisms go in
  Floer (co)homology. One should note that this 
definition  produces an invariant of the
whole manifold and that there are no more open sets involved. 

In this new setting, a
 periodic Hamiltonian $H: \Ss^1 \times M \longrightarrow \RR$ 
is admissible if it satisfies the following: 
\begin{itemize} 
 \item $H \le 0$;
 \item $H \equiv 0$ in a neighbourhood of $\partial M$;
 \item all periodic orbits $x$ of $X_H$ satisfying $\int_x H < 0$ are
 nondegenerate.
\end{itemize} 

There is one important remark that has to be made concerning this
setting: by definition of the admissible Hamiltonians, they have {\it
  plenty} of degenerate $1$-periodic orbits, namely constants in a
neighbourhood of $\partial M$ with action equal to $0$. This has as 
consequence the fact that the truncated homology groups
$FH^*_{]a,\,b]}$ cannot possibly be well 
defined if $0 \in ]a, \, b]$. On the other hand, 
by using admissible almost complex structures  of the type 
described above  they are  well defined if  $a<b<0$ or $0\le a<b$. 
The maximum principle applies in this context  in a direct way to the
parameterized Floer trajectories (which are {\it holomorphic}) near the
boundary and Lemma \ref{max princ}   directly gives the fact  that
(parameterized) 
trajectories stay at a bounded distance from $\partial M$. 

As we now want to let the admissible Hamiltonians go to $-\infty$, 
the partial order will be  
$$H \prec K \qquad \tx{iff} \qquad H(t, \, x) \ge K(t, \, y) \ .$$
If we stick to the cohomological formalism, the inverse  limit
that we have previously encountered
 will therefore change to a direct limit.  
Note the difference with \cite{CFHW} where a slight error concerning
this point has found its way in the paper, with no consequences
whatsoever.

\medskip 
  
\subsubsection
{\sf \large Floer homology of manifolds with contact type boundary.} 
     \label{Homologie symplectique Viterbo} 
     This approach has been developed by 
     C. Viterbo in \cite{functors1}. In comparison with the previous
      ones, it directly isolates the specific features of Floer
     (co)homology that are linked to the Weinstein conjecture by
     building a complex in which to every characteristic corresponds a
     {\it single} (transversally nondegenerate) 
     generator. We define the symplectic completion of $M$
     as 
\begin{equation*}  
\widehat{M} = M \cup_{\partial M \times \{ 1 \} }  
     \partial M \times [1,+\infty[ \ , 
\end{equation*}  
\begin{equation*}  
\widehat{\omega}  =  \left\{  
     \begin{array}{cl}  \omega & \tx{ on  
     } M \\  d(S\lambda\vert) & \tx{ on } \partial M \times [1,+\infty[ \ 
     .  
     \end{array} \right.   
\end{equation*}  
    The glueing is realized through the diffeomorphism 
    (\ref{diffeo Liouville}) induced by the Liouville flow. Formally,
    the admissible Hamiltonians are of the type  $K_\mu$, 
    $\mu > 0$ such that  $K_\mu \equiv 0$ on $M$,  
    $K_\mu   
     (t,p,S) = k_\mu(S) $ on $\partial M \times  [1, + \infty )$,   
     with $k'_\mu(S) =\mu$ for $S \ge 1^+$ and $k_\mu$ 
    convex. We have sketched a typical graph in Figure
    \ref{les 
    trois types de familles admissibles} (3). The homotopies are 
    chosen to be increasing in the $s$ variable and convex. 
    We also note at this point that it would be enough to choose
    Hamiltonians that are linear in the $S$ variable {\it outside a compact
    set} or {\it linear at infinity}. 
 
    The $1$-periodic orbits of  $K_\mu$ are  the constants in the interior
    of  $M$ and the closed characteristics on  $\partial M$ 
    having an action at most equal to  $\mu$. As  $K_\mu$ is time
    independent, any nonconstant orbit is at most transversally 
    nondegenerate and one will have to think in practice to a perturbed
    Hamiltonian or family of Hamiltonians $\widetilde{K}_\mu$. Full
    details concerning the perturbations are provided in  D. Hermann's
    thesis  
    \cite{these Hermann} and the technique of perturbation is due to  
    \cite{FHW}: a nonconstant and transversally nondegenerate
    periodic orbit is to be seen as a Morse-Bott critical 
    nondegenerate circle and a small time-dependent perturbation
    supported in a neighbourhood  will ``break'' the circle  into two 
    nondegenerate periodic orbits, corresponding to the critical points
    of a Morse function on  $\Ss^1$. More
    explanations on this point of view are given section \ref{coh
    Floer boule}, where we perform the same kind of perturbation for
    Morse-Bott nondegenerate spheres, not only circles. 
    The perturbed
    Hamiltonians  $\widetilde{K}_\mu$ satisfy
\begin{itemize}  
  \item $\widetilde{K}_\mu \le \widetilde{K}_{\mu'}$ if $\mu \le 
  \mu'$; 
  \item $(\widetilde{K}_\mu)_{\mu >0}$ is a cofinal family among the
  Hamiltonians   $H: \Ss^1 \times \widehat{M} 
  \longrightarrow \RR$ which verify  $H < 0$ on  $ M 
  \setminus \partial M$.   
\end{itemize}  
     
    The admissible almost complex structures are chosen to be of the
    form (\ref{std a c str}) out of a compact set.
     The existence of  $C^0$ bounds comes from the fact that, 
for generic values of  $\mu$, there are no characteristics of period
    $\mu$ on $\partial M$. As a consequence,  
the  $1$-periodic orbits of  $K_\mu$ are all located in a neighbourhood
of $M$. The Floer trajectories that connect such $1$-periodic orbits
could leave the respective neighbourhood only by having
an interior tangency with a certain  $\partial M \times \{ 
S_0 \}$, a phenomenon which is forbidden by the Lemma below, 
    based again on the maximum principle.

\begin{lem} \label{principe max general}  
{\sf \cite{functors1}} Any solution $u:(D^2(0, \, 1), \, i) 
    \longrightarrow \partial M \times [1-\epsilon, \, \infty [$ of  
     Floer's equation 
     $$u_s + Ju_t -\nabla H(s,\, t, \, u(s, 
    \, t)) = 0$$ 
    with $J$ a standard almost complex structure, 
    $H(s,\, t, \, p, \, S) = h(s, \, t, \, S)$ and  
    $\frac{\partial ^2 h}{\partial s \partial S} \ge 0$ cannot   
    have a strict interior tangency with  some $\partial M \times \{
    S_0 \}$.   
\end{lem}  
 
\demo In what follows we shall denote by $h'$, $h''$ the derivatives
of $h$ with respect to $S$. We put $f=S\circ u$ and we show that it
cannot have an interior strict maximum point. One first computes 
$$\nabla H(s, \, t, \, u) = Cfh'(s, \, t, \, f)\DP{}{S} \ ,$$ 
in view of $\big| \DP{}{S} \big| ^2 = \frac 1 {CS}$. We therefore have
$u_s+Ju_t-Cfh'\DP{}{S}=0$ which, by applying $dS$, gives
$$\partial _s f - C(S\circ u) \lambda(u_t) -Cfh'=0$$
and 
$$\partial _t f +C(S\circ u) \lambda(u_s) =0 \ .$$
We apply $\DP{}{s}$ to the first equation, $\DP{}{t}$ to the second 
and sum up:
\begin{eqnarray*}
  0 & = & \Delta f - C \Big( \partial _s \big( (S\circ u)
  \lambda_u(u_t) \big)  -
  \partial _t  \big( (S\circ u) \lambda_u(u_s) \big) \Big) 
  - C\partial _s \big(fh' \big) \\ 
  & = & \Delta f - C \Big(  u_s\big( (S\lambda)(u_t) \big) - u_t \big(
  (S\lambda) (u_s) \big) \Big) - C\partial _s \big(fh' \big) \\
  & = & \Delta f 
  -C \big(d(S\lambda)(u_s, \, u_t) - (S\lambda)([u_s, \, u_t]) \big) 
- C\partial _s \big(fh' \big) \\
  & = & \Delta f - C \big( |u_s | ^2  - \omega(u_s, \, h'X_{\tx{Reeb}}) \big)
  - C\partial _s \big(fh' \big) \\
  & = & \Delta f - C |u_s | ^2 + C dS(u_s)h' - C \partial _s f h'
  - C f \DP{^2h}{s\partial S} - Cfh''\partial _s f \\
  & = & \Delta f - C |u_s | ^2 -Cf \DP{^2h}{s\partial S}-
  Cfh''\partial _s f \ . 
\end{eqnarray*}

We have used $[u_s, \, u_t]=0$ and $u_t=Ju_s - h'X_{\tx{Reeb}}$ in the
fourth equality. We finally get 

$$\Delta f - Cf
h''(s, \, t, \, f) \partial _s f = C|u_s |
^2 + Cf \DP{^2h}{s\partial S} \ge 0 \ ,$$
with the last inequality holding precisely due to the hypothesis
$\DP{^2h}{s\partial S} \ge 0$. Thus $f$ satisfies the elliptic
inequation of second order without zero order term
$$\Delta f  - C f
h''(s, \, t, \, f) \partial _s f \ge0 \ ,$$
which obviously implies that $f$ cannot have an interior strict
maximum. 

\hfill{$\square$}

\medskip 

We  define as above the groups  $FH^*_{]a, \,
  +\infty[}(M)$, $-\infty \le a < \infty$. 
The new phenomenon is that they are independent of $a$ as soon as
  $a<0$. The explanation is very simple: the terms of the cofinal
  family that we consider all have $1$-periodic orbits with action
  bigger than  $-\delta$, with
  $\delta>0$ arbitrarily fixed.  Indeed, the orbits are of two types:
  on the one hand the critical points of $H$ in the interior of $M$, with
  negative action close to zero, on the other hand the orbits
  corresponding to closed characteristics on 
  $\partial M$, whose action is approximately equal to the (positive)
  area of the latter ones. One should bear in mind at this point that
  the action of a $1$-periodic orbit appearing on a level $S$ is 
  $$A = Sh'(S) - h(S) \ .$$ 
  We therefore have
  $$FH^*(M) = FH^*_{]a, \, +\infty[}(M), \qquad a<0 \ .$$
 
\medskip 

\section{Comments and further properties}    
\label{commentaires}  
  
\bigskip 

\subsection{$C^0$ bounds.} \label{C 0 bounds} 
As we have pointed out repeatedly, the crucial
ingredient in the construction of a Floer (co)homology theory for manifolds
with boundary is the existence of a priori $C^0$ bounds for  Floer
trajectories. The geometric pseudo-convexity 
argument of \cite{functors1} 
has to be modified in a nontrivial manner if one wishes to extend the
class of admissible Hamiltonians and deformations by allowing a more
general dependance on $s$ or $p$ in 
$H(s, \, t, \, p, \, S)$. This motivates the additional work in 
\cite{FH} and  \cite{CFH}, on which we have not given full details. 
Nevertheless, some variant of the maximum
principle enters in both contexts: \cite{FH} Prop. 8   
and \cite{CFH} p. 110. We mention that the ideas in \cite{FH}, pp. 48-56 
can be adapted to the setting of \cite{functors1} in order to establish the
$C^0$ bounds for ``asymptotically linear Hamiltonians'' 
(see \cite{teza mea} \S1.2).

\medskip 
 
\subsection{Relations between  different symplectic homologies.}
\label{relations between homologies} 
\begin{figure}[h]  
      \begin{center}   
      \includegraphics{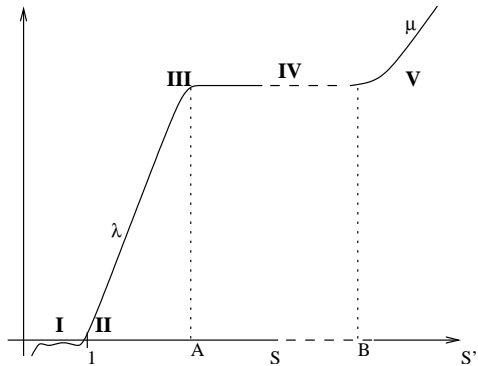}  
      \end{center}   
      \caption{Cofinal family \cite{He} for the isomorphism
        \cite{FH}-\cite{functors1}; $S$ and $S'$ are the coordinates along
        the Liouville field of $U$ and along the radial field in
        $\CC^n$ \label{Ham Hermann}}    
\end{figure}

D. Hermann \cite{He} \S4.3 proves that the homology groups defined in  
    \cite{FH} for a bounded open set  $U \subset  
   \CC^n$  coincide with the ones defined in 
   \cite{functors1} if $U$ has restricted contact type boundary. 
   Moreover, this holds for any fixed range of the action. 
   The proof goes by exhibiting a cofinal sequence of 
   Hamiltonians (Figure \ref{Ham Hermann}) 
   that are admissible in the sense of  \cite{FH}  
   and whose periodic orbits fall in two classes: on one side orbits
   having positive or close to zero action (of type {\bf I} and {\bf
     II}), on the other side orbits
   whose action tends uniformly to  $-\infty$ and which  therefore
   do not count in the Floer complex (of type {\bf III, IV} and {\bf
     V}). This is possible with a suitable choice of the parameters
   $A$, $B$, $\lambda$ and $\mu$. The cofinal sequence is constructed
   so that the orbits of the first kind are also orbits of a
   Hamiltonian that is admissible in the sense of
   \cite{functors1} and the isomorphism between the Floer homology groups
   is induced by an  isomorphism that holds at the level of
   complexes. One consequence of this identification is the
   independence of  $FH^*_{]a, \, +\infty[}(U)$ with respect to  $a<0$
   when $U$ has restricted contact type boundary.
 
I know of no similar construction which identifies the
homology groups defined  in  \cite{CFH} and  \cite{functors1}.  
The  $1$-periodic orbits can no longer be separated in a suitable fashion
by the action as in  \cite{He} and,  
for the same truncation by the action, the complex defined in 
\cite{CFH} generally involves strictly more orbits than the one in 
\cite{functors1}. The presence of extra generators
in the Floer complex raises the probability that ``interesting''
(i.e. geometric) generators may be killed in homology. In this sense 
the homology defined in  \cite{functors1}  
has a more geometric flavour than the one of \cite{CFH}.

One crucial feature of
the definitions in \cite{CFH, CFHW} is that they are intrinsic and
depend only
on the interior of the manifold. This will be crucial in proving the
stability of the action spectrum (see \S\ref{stab}). In comparison, 
the definition in \cite{FH} is inherently extrinsic, while the one in
\cite{functors1} is intrinsic but depends on the boundary as well. 
As a drawback, in the homology of \cite{CFH, CFHW} every geometric 
generator on the boundary appears twice and this may cause some loss
of information when computing the global homology. 

\medskip 

Let us remark at this point that it
would be interesting to find a version of homology
that would depend only on the interior of the manifold and that 
would use a complex in which every geometric generator enters exactly
once. This could lead to topological obstructions for symplectic
fillings, as well as to connections with contact homology.

\medskip 
 
\subsection{Invariance through isotopies.}  Floer homology 
satisfies two important invariance properties, with respect to
deformations of the Hamiltonians and with respect to deformations of
the symplectic forms. We shall state them precisely below. The
Hamitonian invariance is valid under more general hypothesis than the
ones adopted here, in particular for some 
homotopies that are not necessarily increasing. 
But for our needs the results below will do. Note that it
does not matter whether one works with the definition in
\cite{CFH} or with that in \cite{functors1}. 
\begin{thm} {\sf \cite{functors1}} \label{premiere invariance par isotopie} 
Let $H(s)$, $s\in [0, \, 1]$ be an {\rm admissible deformation} 
(cf. \S\ref{hom Floer type contact}) such that the $1$-periodic orbits
of $X_{H(s)}$, $s\in [0, \, 1]$ all stay in a fixed compact set. 
Let $a_s$, 
$b_s$ be continuous families of parameters in $\overline{\RR}$ such
  that  $H(s)$ has no $1$-periodic orbits with action equal to  
$a_s$ or $b_s$. There is a canonical isomorphism 
$$FH^*_{]a_1, \, b_1]}(H_1) \stackrel{\sim}\longrightarrow FH^*_{]a_0, 
\, b_0]}(H_0) \ ,$$ 
obtained in the usual manner by considering solutions of the equations 
 (\ref{Fl param 1} - \ref{Fl param 2}). 
\end{thm}  
One should note that we suppose in particular 
$H_0 \prec H_1$. The requirement that all $1$-periodic orbits stay in
a fixed compact set is essential. We leave it to the reader to extract
a ready-made example from the
computations given in the next section. 

The invariance of Floer homology with respect to deformations of the
symplectic form is closely related to the preceding result. It can be
given a heuristic interpretation  by taking  $a_s 
\equiv -\infty$, $b_s \equiv +\infty$:  
\begin{thm} {\sf \cite{functors1}} \label{deuxieme invariance par isotopies} 
   Let $M$ be a compact symplectic manifold with boundary and  
   $\omega_s$, $s \in [0, \, 1]$ an isotopy of symplectic forms with
   respect to which the boundary  $\partial M$ is of contact
   type. The nontruncated Floer (co)homology groups corresponding to
   $\omega _0$ and $\omega_1$ are isomorphic: 
   $$FH^*(M, \, \omega_0) \simeq 
   FH^*(M, \, \omega_1) \ .$$ 
\end{thm} 
    The isomorphism is again obtained by counting solutions of an
    equation that is very much similar to (\ref{Fl param
    1}). It is crucial to work with the nontruncated homology in
    order for the morphism to be well defined. It is also important to
    note that the isotopy must not necessarily preserve the cohomology
    class of the symplectic form. Let us also bear in mind the
    following useful corollary which shows that the nontruncated
    Floer homology of a manifold with boundary $M$ is in fact an
    invariant of its completion 
$\widehat{M}$. The result is stated for open sets with smooth
    boundary only for the sake of applying the above theorem as
    such, although it is true for general open sets if one
    suitably defines the Floer homology groups. 
    As a side remark, note that the invariance with respect
    to deformations of the symplectic form raises the question of the
    {\it differentiable} invariance of Floer homology of manifolds
    with boundary. 
\begin{cor} \label{corollaire gi} 
    Let  $U, \ U'$ be two open sets  {\rm with smooth boundary}  
    in $\widehat{M}$ which satisfy the following conditions: 
    \begin{itemize}  
     \item $\partial U$ and $\partial U'$ are contained in the domain
       of definition of the Liouville vector field extended to
       $\widehat{M}$   
       as $S\DP{}{S}$;  
     \item the Liouville vector field is transverse and outward
       pointing along   
       $\partial U, \ \partial U'$.   
    \end{itemize}  
    Then    
    $$ FH^*(U) \simeq FH^*(U') \ .$$  
\end{cor}  
   \demo     
   One can realise a differentiable isotopy between  $U$ and $U'$  
   along the Liouville vector field. This corresponds to an isotopy of
   symplectic forms on  $U$ starting from its initial   
   symplectic form $\omega_0$ and ending to the one induced from 
   $U'$, denoted by $\omega_1$.  The invariance theorem ensures 
   $FH^*(U,\, \omega_0) \simeq  
   FH^*(U,\, \omega_1)$. But $(U, \, \omega_1)$ and $(U', \,
   \omega_0)$ are symplectomorphic and hence  we also have $FH^*(U,  
   \, \omega_1) \simeq FH^*(U', \, \omega_0)$. 
  
   \hfill{$\square$}

\section{A computation: balls in $\CC^n$}

\label{coh Floer boule} The Floer cohomology groups   
$FH^*_{]a, \, b]}$ of the ball
 $D^{2n} \subset \CC^n$ have been computed for the first time in 
 \cite{FHW},  
together with those of ellipsoids and polydiscs, for arbitrary 
 $a, \, b \in 
\overline{\RR}$. The nontruncated cohomology groups vanish:
 $FH^*(D^{2n}) = 0$. The method of 
 \cite{FHW} consists in approximating  the ball 
 by ellipsoids whose characteristics are transversally 
 nondegenerate. The homology of the latter is studied by considering
 perturbed Hamiltonians whose nontrivial 
 $1$-periodic orbits are nondegenerate and come in pairs which
 correspond to 
 transversally nondegenerate characteristics. The two $1$-periodic
 orbits in a pair correspond to the two critical points of a Morse
 function on the corresponding 
 characteristic, which is identified with $\Ss^1$. 
 The difference between their Conley-Zehnder
 indices is equal to 
$1 = \dim \Ss^1$. The same perturbation technique is used in 
 \cite{CFHW} in the proof of  the stability of the action spectrum
 of the boundary of a symplectic manifold.

The closed characteristics on the sphere  
$\Ss^{2n-1} = \partial D^{2n}$ are the great circles together with their
positive multiples. They appear in families that are parameterized by
the sphere  $\Ss^{2n-1}$ itself. We shall present in section \ref{ball
  by perturbing} a way of perturbing the natural Hamiltonians
$H(z)=h(|z|^2)$   
such that every  manifold of
closed characteristics produces two 
nondegenerate periodic orbits, whose difference of
Conley-Zehnder indices is equal to   
$2n-1=\dim \Ss^{2n-1}$. This is the higher dimensional 
analogue of the perturbation
method in \cite{FHW}. The interest of this point of view lies in
that it allows one to understand the vanishing of the Floer
homology of the ball without using as an intermediate ingredient
the homology of ellipsoids. It also involves some interesting
computations of Conley-Zehnder indices and may serve as a hands-on
example for the constructions in \cite{Po} \S3.4 and 
\cite{BPS} \S5.2 which compute the local Floer homology of a
Lagrangian intersection and, respectively, of a Hamiltonian admitting
a Morse-Bott nondegenerate 
manifold of $1$-periodic orbits.

\subsection{Recollections on the Robbin-Salamon index} 
We remind that the Conley-Zehnder index is an integer that is
associated to a path of symplectic matrices starting at $\tx{Id}$ and
ending at a matrix $\Psi$ such that $\det(\tx{Id}-\Psi) \neq 0$. 
This is precisely the situation  encountered after linearizing the
Hamiltonian flow along a nondegenerate periodic orbit. 

We shall compute below some Conley-Zehnder indices through the
intermediate of the Robbin-Salamon index \cite{RS1}. The latter 
generalizes the Conley-Zehnder index to arbitrary paths 
  in  $\tx{Sp}(2n, \,  \RR)$. The reader can consult 
\cite{Sa} \S2.4. for a 
  review of the basic properties of the Conley-Zehnder index and for
  its computation using crossing numbers, as well as  \cite{RS1} for
  a complete account on the Robbin-Salamon index. We shall sketch here
  only the facts that are used in the sequel. 

 The Robbin-Salamon index $i_{RS}(\Psi)$ 
 of an {\it arbitrary} path of symplectic matrices $\Psi: [0,
 \, 1] \longrightarrow \tx{Sp}(2n)$ is defined as the index of the
 Lagrangian path $\tx{gr}(\Psi)$ in $(\RR^{2n}\times \RR^{2n}, \,
 (-\omega) \oplus \omega)$ {\it relative to the diagonal} $\Delta =
 \{(X, \, X) \ : \ X \in \RR^{2n} \}$. Its computation makes use of
 the notions of {\it simple crossing} and {\it crossing form}. 

 A {\it crossing} is a number $t\in [0, \, 1]$ such that
 $\tx{gr} \, \Psi(t) \cap \Delta \neq 0$ or else $\ker(\tx{Id} -\Psi(t))
 \neq 0$. The {\it crossing form} at a crossing $t_0$ is the quadratic
 form $\Gamma(\Psi, \, t_0): \ker(\tx{Id} -\Psi(t))
 \longrightarrow  \RR$ defined by
 $$\Gamma(\Psi, \, t_0)(v) = \omega_0(v, \, \dot\Psi(t_0)v) \ .$$
 Any path in $\tx{Sp}(2n, \, \RR) $ is solution to a differential
 equation $\dot\Psi(t) = J_0 S(t) \Psi(t)$, with $S(t)$ a symmetric
 matrix. One can therefore write 
 $$\Gamma(\Psi, \, t_0)(v) = \langle v, S(t_0) v \rangle \ . $$ 
 A crossing $t_0$ is called {\it simple} if the crossing form
 $\Gamma(\Psi, \, t_0)$ is nondegenerate. A simple crossing is
 necessarily isolated. 

 Any path $\Psi$ is homotopic with fixed end points to a path
 having only simple crossings. The index $i_{RS}(\Psi)$ of a path
 having only simple crossings is defined to be 
 $$i_{RS}(\Psi) = \frac 1 2 \, \tx{sign} \, \Gamma(\Psi, \, 0) + 
\sum _{0 < t < 1 } \tx{sign} \, \Gamma (\Psi, \, t) +  
\frac 1 2 \, \tx{sign} \, \Gamma(\Psi, \, 1) \ , $$
 where the summation runs over all simple crossings of $\Psi$. We
 recall that the signature of a nondegenerate quadratic form is
 the difference between the number of its positive eigenvalues
 and the number of its negative eigenvalues. 

 The main features of the Robbin-Salamon index are the following. 
\renewcommand{\theenumi}{\roman{enumi}}
\begin{enumerate}
\item {\it additivity under concatenations of paths:}
  $i_{RS}(\Psi|_{[a, \, b]}) + i_{RS}(\Psi|_{[b, \, c]}) =
  i_{RS}(\Psi|_{[a, \, c]}) $ ; 
\item $i_{RS}$ characterizes paths up to homotopy with fixed end
  points ;
\item {\it additivity under products:} 
 $i_{RS}(\Psi ' \oplus \Psi '') = i_{RS}(\Psi ') + i_{RS}(\Psi'')$ .
\end{enumerate}

\subsection{Direct computation of the cohomology of a ball.} The
argument that we present below appears in \cite{FHW} and
\cite{functors1}.  
Let us identify $\CC^n$ with $\RR^{2n}$ by associating to  $z =
x+ iy$, $x, \, y \in \RR^n$ the vector  $(x_1, \ldots, \, x_n, \, y_1,
\ldots, \, y_n)$.  
We consider on  $\CC^n$ the standard symplectic form  $\omega_0 =
\sum _{i=1} ^n dx_i \wedge dy_i$ with its primitive $\lambda_0 = \frac
1 2 \sum_{i=1}^n x_i dy_i - y_i dx_i$. The associated Liouville
vector field  $X(z)= \frac 1 2 z$, $z\in \CC^n$ is transverse to all
the spheres that are centered at the origin. We shall use the standard
almost complex structure given by complex multiplication with $i$ i.e.
$J(x_1, \ldots, \, x_n, \, y_1,
\ldots, \, y_n) = (-y_1,
\ldots, \, -y_n, \, x_1, \ldots, \, x_n)$. 
The Reeb vector field on  $\Ss^{2n-1}$ is  
$X_{\tx{Reeb}}(z)=2Jz$ and its closed orbits are the great circles and
their positive multiples, with action  $k\pi$, $k \in \mathbb{N}^*$. 
It is clearly possible to construct a cofinal family of Hamiltonians
of the form  $H_\lambda(z) = \rho(|z|^2)$, with   
$\rho:[0, \, \infty[ \longrightarrow \RR$ satisfying  
$\rho |_{[0, \, 1]} <0$, $\rho'' \ge 0$, $\rho' | 
_{[1+ \epsilon, \, \infty[} = \lambda$, $\rho' | _{[0, \, 
1+\epsilon)} >0$, $\epsilon = \epsilon(\lambda) >0$. When $\lambda 
>0$ is not a multiple of $\pi$ the  $1$-periodic orbits of  
$H_\lambda$ are parameterized by the critical point  $0$ (the constant
orbit) and the spheres 
 $\{z \, : \, \rho'(|z|^2) = l\pi \}$, $1\le l \le k$ 
where $k\pi < \lambda < (k+1)\pi$.
 
Let therefore  
$\lambda \neq k\pi$ and take $b$ sufficiently large compared to  
$\lambda$. We claim that 
\begin{equation} \label{petit iso intermediaire}  
FH^*_{]-\infty, \, b]}(H_\lambda ) = FH^*_{]-\infty, \, 
b]}(\lambda|z|^2 - c) \ , 
\end{equation}  
where $0< c \le 1$ is a constant such that 
$\lambda|z|^2 -c \le H_\lambda$. Note that it is always possible to
choose $c\le 1$ if $\lambda$ is large enough, due to the cofinality of 
$(H_\lambda)_\lambda$.  
We can connect $\lambda|z|^2-c$ with $H_\lambda$ through an
admissible homotopy of Hamiltonians of the same type as
$H_\lambda$. The periodic orbits that are created correspond to closed
characteristics on 
$\Ss^{2n-1}$ and we easily see that their action will not cross $b$ if
the latter is large enough (e.g. $b > \lambda +1$ for our choice of
$c$). 
The isotopy invariance theorem  \ref{premiere 
invariance par isotopie} implies (\ref{petit iso intermediaire}).  
 
In order to fix ideas, let us suppose that  $k\pi < \lambda < (k+1) \pi$.  
The Hamiltonian $Q_\lambda= \lambda|z|^2 - c$ has a unique
$1$-periodic orbit: the constant  $z=0$. We have seen above that the
natural grading on Floer (co)homology is given by the Conley-Zehnder
index taken with opposite sign. This amounts to the computation of the
index of the $1$-periodic orbits of the vector field 
 $-X_{Q_\lambda}(z) =
2\lambda Jz$, whose flow is 
$\varphi_t(z)=e^{2i\lambda t}z$. Its linearization is 
$d\varphi_t (z)\cdot Y = e^{2i\lambda t}Y$, or  
$$d\varphi_t = \left( \begin{array}{ccc}  
\left( \begin{array}{cc} \cos  2\lambda t & - \sin 2\lambda t \\   
\sin 2\lambda t & \cos 2 \lambda t \end{array} \right)   
& & 0 \\ 
& \ddots & \\ 
0 & &  
\left( \begin{array}{cc} \cos  2\lambda t & - \sin 2\lambda t \\   
\sin 2\lambda t & \cos 2 \lambda t \end{array} \right)   
\end{array}  
\right)$$ 
The Conley-Zehnder index is 
$i_{CZ}(d\varphi_t) = n \cdot  i_{CZ}\Big( \left( \begin{array}{cc} \cos 
2\lambda t & - \sin 2\lambda t \\   
\sin 2\lambda t & \cos 2 \lambda t \end{array} \right)  \Big) = n 
(2k+1)$~: the path $t \mapsto e^{2i\lambda t}$ in $\tx{Sp}(2, \, 
\RR)$ has $k$ interior crossings at  $t=l\pi/\lambda$, $1\le l \le k$ 
and one initial crossing at  $t=0$, with an intersection form of
constant signature equal to  $2$. We get thus the nontruncated
cohomology 
$$FH^*(H_\lambda)= \left\{ \begin{array}{rl} \ZZ, & * = n(2k+1) \, ,\\ 
0, & *\neq n(2k+1) \end{array} \right.  $$  
 
As a consequence we have  
$$\displaystyle FH^*(D^{2n}) = \lim_{\longleftarrow}
\, FH^*(H_\lambda) = 0 \ .$$

 \begin{figure}[h]  
      \begin{center}   
      \includegraphics{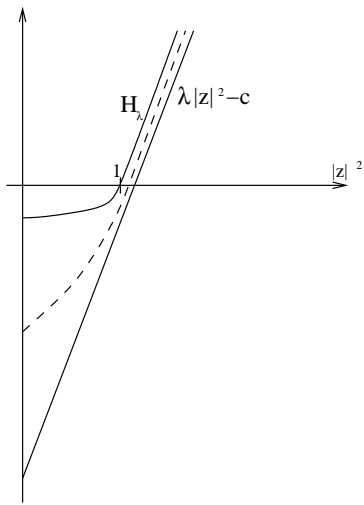}  
      \end{center}   
      \caption{Computation of the Floer homology of
        $D^{2n}$\label{figure calcul boules}}    
\end{figure}    
 
\subsection{Alternative point of view: perturbation of critical
  Morse-Bott manifolds.} \label{ball by perturbing}
We present now an alternative point of view which consists in
  directly perturbing the critical manifolds of  $H_\lambda$. Our main
  computation concerns the (opposite) 
  Conley-Zehnder indices of the corresponding
  generators in the Floer complex.

Let us first compute the  Robbin-Salamon indices $i_{RS}$ for the 
$1$-periodic orbits of  $H_\lambda$, $k\pi < \lambda < 
(k+1)\pi$ considered with the opposite orientation. 
We have  $- X_{H_\lambda}(z)=2\rho'(|z|^2)Jz$ and the flow is  
$\varphi_t(z)=e^{2\rho'(|z|^2) Jt}z$. The linearization of the flow is   
$$d\varphi_t(z) \cdot Y = e^{2\rho'(|z|^2) Jt}\cdot Y + 
2\rho''(|z|^2)Jt\cdot e^{2\rho'(|z|^2) Jt}\cdot 2\langle z, \, Y 
\rangle \cdot z$$ 
For $z=0$ the linearization is  $d\varphi_t(0) = e^{2\rho'(0)Jt}$ and  
satisfies the differential equation  $\dot{\Psi}(t) = J 
\, \tx{diag}(2\rho'(0)) \, \Psi(t)$, $\Psi(0)= \tx{Id}$. For a
sufficiently large $\lambda$ 
the cofinality of  $H_\lambda$ imposes 
$\rho'(0)<\pi$. As a consequence the path $e^{2\rho'(0)Jt}$ has a
unique crossing at  $t=0$ with intersection form of signature
$2n$. We get   
\begin{equation}  
i_{RS}(z\equiv 0)=n \ . 
\end{equation}  
Let us now look at a sphere  $S_l=\{ z \, : \, \rho'(|z|^2) =
l\pi \} $, $1  
\le l \le k$. Consider the $J$-invariant symplectic decomposition  
of $T_z\RR^{2n}$, $z\in S_l$ given by   
$$T_z\RR^{2n} = \RR \langle \frac 1 {2 l\pi |z|^2} z \rangle 
\ \oplus \ \RR \langle 2 l\pi Jz \rangle \ \oplus \ \xi_z \ ,$$ 
where $\xi=TS_l \cap J TS_l$ 
is the contact distribution on  $S_l$. The matrix of  
$d\varphi_t$ is written with respect to this decomposition as 
$$d\varphi_t (z) = \left(  
\begin{array}{cc}  
\left( \begin{array}{cc} 1 & 0 \\ 
\frac {\rho''} { l^2\pi^2 } \cdot t & 1 \end{array} \right)  
&  \\ 
& d\varphi_t |_{\xi_z}   
\end{array}  
\right) 
$$ 
As $d\varphi_t(z)|_{T_z S_l} = e^{2l\pi Jt}$ we can still write in the
canonical trivialization of  $T\RR^{2n}$:  
$$d\varphi_t(z) = \chi(t)\circ \Psi(t) \ ,$$ 
with $\Psi(t) = e^{2l\pi Jt}$ and
$$\chi(t) = \left(  
\begin{array}{cc}  
\left( \begin{array}{cc} 1 & 0 \\ 
\frac {\rho''} {l^2\pi^2} \cdot t & 1 \end{array} \right)  
&  \\ 
& \tx{Id}_{2n-2} 
\end{array}  
\right) 
$$ 
The homotopy (cf. \cite{CFHW}) 
$$K(s, \, t) = \left \{ \begin{array}{ll} \chi(st) \Psi(\frac {2t} {s+1}) 
, & t \le \frac {s+1} 2 \\ 
\chi\big( (s+2)t - (s+1) \big) \Psi(1) , & \frac{s+1} 2 \le t  
\end{array} \right. 
$$ 
connects with fixed ends the path $\chi(t)\Psi(t) $ to the concatenation
of 
$\Psi(t) $ and $\chi(t)\Psi(1)$. We thus get  
$$i_{RS}(d\varphi_t(z)) =  i_{RS}(\Psi(t)) + i_{RS}(\chi(t)\Psi(1))  = 
2ln + \frac 12  \ ,$$ 
or else  
\begin{equation}  
i_{RS}(S_l) = 2ln + \frac 12 \ . 
\end{equation}  
 
Once we have computed the indices before perturbing, we can go on with
the description of the latter. By a {\it time-dependent} change of
variables we can suppose that 
$$X_{H_t} \equiv 0 \quad \tx{on } \quad S_l \ .$$ 
We leave it to the reader to work it out or check
\cite{CFHW}, p.34. The idea is to spin $\Ss^l$ the other way round so
that all its points become fixed under the new flow.
 
Let us now choose a Morse function having exactly two critical points 
$h: S_l \longrightarrow \RR$. Let  $H_\delta = H+ \delta 
h$, $\delta > 0$. One can show by an argument similar to 
\cite{CFHW} that, for $\delta$ small enough,   
$X_{H_\delta}$ has precisely two $1$-periodic (constant) orbits in the
neighbourhood of  $S_l$, that correspond to the two critical points of 
$h$. Let  $\Psi(t)$ be the linearized flow of  $-X_H$, $\Phi_\delta(t)$ 
be the linearized flow of  $-X_{\delta h}$, $\Psi_\delta(t)$ be the
linearized flow of  $-X_{H_\delta}$. For small $\delta$ the paths  
$\Psi_\delta(t)$ and $\Psi(t)\Phi_\delta(t)$ are homotopic with ends
in  $\tx{Sp}_0(2n, \, \RR)$, the set of symplectic matrices with
eigenvalues different from $1$. The homotopy is given by 
$L(s, \, t) = 
\Psi_{s\delta}(t)\Phi_{(1-s)\delta}(t)$. The same argument as above
shows that 
$$i_{RS}(\Psi_\delta) = i_{RS}(\Psi) + i_{RS}(\Psi(1)\Phi_\delta) \ 
.$$ 
But $i_{RS}(\Psi)=\frac 1 2 +2ln$ by the above computations. It is
therefore enough to find the index of 
$M(t)=\Psi(1)\Phi_\delta(t)$ and we again use the definition.  
 We can extend  $h$ to 
$\widetilde{h}$ in the neighbourhood of  $S_l$ by  
$\widetilde{h}(z, \, S) = h(z)$, where  $(z, \, S) \mapsto Sz$ is a
parameterization of a neighbourhood of  $S_l$ 
by $S_l \times [1-\epsilon,
\, 1+\epsilon]$.  
Then $-X_{\delta h}(z) = \delta J \nabla \widetilde{h}(z) = \delta J 
\nabla h(z)$, $z\in S_l$ and the linearization of the flow at the
critical point 
$z_0 \in S_l$ is a solution of the differential equation 
$\dot{A}(t)=D\big(-X_{\delta h}(z_0)\big) \cdot A(t)$, $A(0) = 
\tx{Id}_{2n}$. This implies
$$d \varphi_t^{-X_{\delta h}}(z_0) = e^{\delta t J \nabla ^2 
\widetilde{h}(z_0)} \ . $$ 
But it is easy to see that  $\nabla ^2 \widetilde{h} (z_0) = \left(  
\begin{array}{cc} \nabla ^2 h(z_0) & 0 \\ 0 & 0 \end{array} \right)$ 
in the decomposition $T_{z_0}\RR^{2n} = T_{z_0}S_l \, \oplus \, \RR  
\langle z_0 \rangle $. This immediately implies 
$\ker d\varphi_t^{-X_{\delta h}} = \{ 0 \} $ for $t >0$ and one is reduced to
compute the intersection form at  $t=0$. We have   
$$ \ker(\tx{Id} - \Phi_{\delta}(0) =  T_{z_0}S_l  $$
and therefore
$$\tx{sign} \, \Gamma(\Phi_{\delta}, \, 0) = \tx{sign} \, \delta
\nabla ^2 h(z_0) \ .$$
The signature of  $\Gamma$ is therefore equal to 
$2n-1$ at the critical point of index  $0$, respectively equal to 
$-(2n-1)$ at the critical point of index  $2n-1$. This means  
$i_{RS}(\Psi(1)\Phi_\delta(t)) = \pm (n-\frac 1 2 )$ and
$$i_{RS}(d\varphi_t^{-X_{H_\delta}}(z_0)) = \left\{ \begin{array}{ll} 
2ln +n , & \tx{ind}_{\tx{Morse}}(z_0) = 0 \\ 
2ln -n +1, & \tx{ind}_{\tx{Morse}}(z_0) = 2n-1 \end{array} \right. $$ 
 
The situation is therefore the following: for
$k\pi < \lambda < (k+1)\pi$ the Hamiltonian  $H_\lambda$ admits a  
constant orbit $z\equiv 0$ of index $n$ and  $k$ manifolds of periodic
orbits diffeomorphic to  $\Ss^{2n-1}$ and having indices  
$2ln +\frac 1 2$, $1 \le l \le k$. After perturbation, each such
manifold produces two nondegenerate orbits of indices 
$2ln - n
+1$ and $2ln +n$, corresponding respectively to the maximum and
minimum of the perturbing function.
The truncated homology  $FH^*_{(-\infty, \, b]}(H_\lambda)$, $b 
> \lambda$ was already computed above and it is nontrivial only in
degree $n+2kn$, where it equals  $\ZZ$. We infer that the arrows in
the perturbed Floer complex run as shown below, where the third line
indicates the degree in the Floer complex, given by the opposite of
the Conley-Zehnder index
 and the first line indicates the action of
the corresponding orbit (compare with \cite{FHW}). 

{\scriptsize
\begin{equation}\label{diag important}
\xymatrix{  
\epsilon & \pi - \epsilon & \pi+\epsilon & 2\pi -\epsilon & &
(k-1)\pi-\epsilon & (k-1)\pi +\epsilon & k\pi -\epsilon \\
\ZZ \ar@{->}[r]^{\tx{ Id}} & \ZZ & \ZZ \ar@{->}[r]^{\tx{ Id}} & \ZZ 
& \ldots & \ZZ \ar@{->}[r]^{\tx{ Id}} & \ZZ & \ZZ \\ 
n & n+1 & 3n & 
3n+1 & & (2k-1)n & (2k-1)n+1 & 
(2k+1)n  
}
\end{equation} 
}
 
This allows one to also compute the values of the truncated Floer
cohomology groups $FH^*_{]a, \, b]}(D^{2n})$ for arbitrary $a, \, b \in 
\overline{\RR}$, recovering  the results in \cite{FHW}.

\medskip

\section{Applications}

As we have already pointed out, each of the three constructions
described above has specific features that allow particular
applications. We shall present below the symplectic 
classification of ellipsoids
and polydiscs using symplectic homology as in \cite{FH}, the stability
of the action spectrum of the contact type boundary of symplectic manifolds
 \cite{CFHW}, as well as 
applications to
Weinstein's conjecture by using the  homology defined in 
\cite{functors1}. These three problems will make the reader familiar
with some techniques that are representative for the field.

\subsection{Polydiscs and ellipsoids} \label{poly ell}

We use in this section the homology groups defined by Floer and Hofer
or those defined by Viterbo (cf. sections \ref{homologie FH},
\ref{Homologie symplectique Viterbo}, \ref{relations between
  homologies}). We  explain the following two theorems. 

\begin{thm} \label{ellipsoids} 
Let $r=(r_1, \, r_2, \, \ldots, \, r_n) \in \RR_+^*$, $r_1
  \le r_2 \le \ldots \le r_n$ and denote 
  $$E(r) = \big\{ (z_1, \ldots, \, z_n) \in \CC^n \ : \ \sum_{i=1}^n
  \frac {|z_i|^2}{r_i^2} <1 \big\} \ .$$
 Then $E(r)$ and $E(r')$ are symplectically diffeomorphic
 if and only if $r=r'$. 
\end{thm} 

\begin{thm} \label{polydiscs}
Let $r=(r_1, \, r_2, \, \ldots, \, r_n) \in \RR_+^*$, $r_1
  \le r_2 \le \ldots \le r_n$ and denote 
$$D^{2n}(r) = B^2(r_1) \times \ldots \times B^2(r_n) \subset \CC^n \
  .$$
Then $D^{2n}(r)$ and $D^{2n}(r')$ are symplectically diffeomorphic if
  and only if $r=r'$. 
\end{thm} 

Both theorems follow from the explicit computation of the {\it
  truncated} symplectic (co)homology of ellipsoids and polydiscs and
  we shall address the case of ellipsoids. The computation for
  polydiscs is of a similar nature and we refer the interested reader
  to  \cite{FHW}, p. 583. 

{\it Sketch of proof for Theorem \ref{ellipsoids}.} One interesting
remark concerning the Floer complex that we have computed in
\S\ref{coh Floer boule} is that the same homology would have been
obtained out of the following complex: 

{\scriptsize
\begin{equation} \label{complexe ellipsoides boule}
\xymatrix{
0 \ar[r] & (\ZZ, \, n; \, 0) \ar[r]^{\tx{Id}} & (\ZZ, \, n+1; \, \pi)
\ar[r]^0 & (\ZZ, \, n+2; \, \pi) \ar[r]^{\tx{Id}} & \ldots \ar[r]^0 & 
(\ZZ, \, 3n; \, \pi) \ar[r]^{\tx{Id}} & (\ZZ, \, 3n+1; \, 2\pi)
\ar[r]^0 & \ldots }
\end{equation} 
}

The meaning of the notation is the following: 
in a term $(\ZZ, \, k; \, \alpha)$ we have $k$ standing for the grading and
$\alpha$ for the corresponding action of the generator. If one wishes
to compute the truncated cohomology $FH^*_{]a, \, b]}$ 
    the complex to consider is $C_b/C_a$, with $C_b =
    \oplus_{\alpha \le b}
    (\ZZ, \, k; \, \alpha)$. This corresponds to formally replacing 
the arbitrarily small  $\epsilon$ in (\ref{diag
  important}) with $0$. 

The complex (\ref{complexe ellipsoides boule}) is precisely 
the one used in \cite{FHW}
and arises in a natural  geometric way. To the difference of our
method of perturbing the spheres of characteristics with action
$k\pi$, $k\in \ZZ_+^*$ and produce $2$ orbits, the method of the
original paper is to approximate the ball by a generic ellipsoid
having $n$ simple characteristics $(0, \ldots, r_je^{2\pi i t/r_j},
\ldots, \, 0)$ whose areas $\{\pi r_j^2 \}$, $1\le j \le n$ are 
linearly independent over $\QQ$, and then perturb {\it each} such
characteristic in order to finally produce $2n$ nondegenerate
orbits. For the perturbations of the simple characteristics that arise
for an approximation of the sphere the indices are precisely $n+1, \, 
n+2, \ldots, \, 3n$. 

This discussion is of course no proof of the general case, but gives a
geometric explanation for the construction of the complex that computes
the symplectic cohomology of a general ellipsoid, where every
closed characteristic (simple or not) gives rise to $2$ orbits having
action arbitrarily close to some $k\pi r_j^2$, $k\in \ZZ_+^*$.   

{\scriptsize
$$FH^*_{]a, \, b]}(E(r)) = H^*(C_b/C_a) \ ,$$
\begin{equation} \label{complexe ellipsoides} 
C_b = \xymatrix{
0 \ar[r] & (\ZZ, \, n) \ar[r]^{\tx{Id \ }} & (\ZZ, \, n+1)
\ar[r]^0 & (\ZZ, \, n+2) \ar[r]^{\tx{\quad Id}} & \ldots \ar[r]^{0 \qquad
  \quad} &  
(\ZZ, \, n+ 2m(b; \, r)) \ar[r] & 0 \ar[r] & 0 \ar[r] & \ldots 
}
\end{equation} 
}

where $$m(b; r)= \# \big\{ (k, \, j) \in \ZZ_+^* \times \{ 1, \ldots, \, n
\} \ : \ k\pi r_j^2 \le b \big\} \ .$$
 
The reader can easily get convinced that, for the ball of radius one,
the two complexes (\ref{complexe ellipsoides boule}) and
(\ref{complexe ellipsoides}) are the same. 

In the case of a generic ellipsoid Theorem \ref{ellipsoids} follows
now at once, as one can easily compute for example 
$FH^{n+2j}_{]\pi r_j^2
-\epsilon, \, \pi r_j^2 +\epsilon] } \simeq \ZZ$ and $FH^*_{]\pi r_j
^2 -\epsilon, \, \pi r_j ^2 +\epsilon]} = 0$ if $*\neq n+2j$. 

\hfill{$\square$}

\subsection{Stability of the action spectrum} \label{stab}

In this section we use the version of symplectic homology described in
\cite{CFHW} (cf. section \ref{CiFH}). 
The ``stability of the action spectrum'' is a problem related to the
question of the extent to which the interior of a symplectic manifold
determines its boundary. One answer is that, in the nondegenerate
case, the set of values of the areas of closed characteristics on the
boundary is determined by the interior of the manifold. We
first introduce the relevant definitions. 

\begin{defi} 
 A symplectic manifold $(M, \, \omega)$ is said to satisfy the
 {\rm symplectic asphericity} condition if $\langle \omega, \,
 \pi_2(M) \rangle = 0$. If in addition  
 $\langle c_1(M), \, \pi_2(M) \rangle =0$ we say that the manifold is
 {\rm strongly symplectically aspherical} (the Chern class is
 computed with respect to an almost complex structure $J$ that is
 compatible  with the symplectic form, and does not depend on $J$).   
\end{defi} 

\begin{defi} 
Let $(M, \, \omega)$ be a
compact symplectic manifold with nonempty contact type boundary.
Assume $(M, \, \omega)$ is  strongly symplectically aspherical. 
We note 
$X$ the Liouville vector field, $\lambda$ the Liouville form and
$X_{\tx{Reeb}}$ the Reeb vector field. 

Let $x: [0, \, T] \longrightarrow \partial M$ 
be a contractible closed characteristic with  $x(0) = x(T)$. 
We define the {\rm action of $x$}  to be 
$$A(x) = \int_0^T x^*\lambda \ .$$

We define the {\rm index $i_{RS}(x)$ of $x$} to be the Robbin-Salamon index
 of the path 
$$\Gamma(t) = \Psi(e^{2\pi i t/T}) \circ 
{\psi_t}_* \in \tx{Symp}(T_{x(0)}M) \ ,$$
where $\psi_t$ is the flow of the Reeb vector field on $\partial M$,
${\psi_t}_*:T_{x(0)}M \longrightarrow T_{x(t)}M$ 
denotes its linearization extended by ${\psi_t}_* (X(x(0))=X(x(t))$,
and $\Psi:TM|_{D^2}\longrightarrow D^2 \times T_{x(0)}M$ 
denotes a trivialization of $TM$ over a filling disc $D^2$ for
$x$. 
\end{defi} 

One should note that the condition $\langle \omega, \, \pi_2(M)
\rangle =0$ ensures that the value of the action is independent on the
Liouville form, as we have $\int_0^Tx^*\lambda =
\int_{D^2}\bar{x}^*\omega$, with $\bar{x}:D^2 \longrightarrow M$ such
that $\bar{x}(e^{2\pi i t/T})=x(t)$. On the other hand, the hypothesis
$\langle c_1, \, \pi_2(M) \rangle = 0$ ensures that any two
trivializations of $TM$ over disks with common boundary are homotopic
along the boundary and this means that the value  of
$i_{RS}(x)$  does not  depend on the choice of trivialization.

\begin{defi}\cite{CFHW} 
Let $(M, \, \omega)$ be a compact symplectic manifold with 
nonempty contact type boundary, satisfying the symplectic asphericity
  condition $\langle \omega, \, \pi _2(M) \rangle = \langle c_1, \,
  \pi _2(M) \rangle =0$. Denote by $\mc C(M)$ the set of contractible 
  closed
  characteristics on $\partial M$. The action spectrum $\mc A(\partial
  M)$ is defined as 
\begin{equation}
  \mc A(\partial M) = \big\{ ( A(x), \, -i_{RS}(x)) \ | \ x\in \mc C(M)
  \big\} \ .
\end{equation} 
\end{defi}

\begin{defi} Let $x: [0, \, T] \longrightarrow \partial M$ be a closed
  characteristic of period $T$, parameterized by the Reeb vector
  field. We say that $x$ is {\rm transversally nondegenerate} if the
  restriction of the linearization of the Reeb flow 
$${\psi_T}_* : \xi_{x(0)} \longrightarrow \xi_{x(0)} $$
 has no eigenvalue equal to $1$.
\end{defi}

We will see below that, under this hypothesis, the index $i_{RS}(x)$
is an integer.

\begin{thm} \cite{CFHW} \label{thm stab}
 Let $M$ and $N$ be compact symplectic manifolds with nonempty
 boundary of contact type, satisfying the symplectic asphericity
 condition and having transversally nondegenerate closed
 characteristics.  Assume that their interiors 
 $\dot{M}$, $\dot{N}$ are  symplectically
 diffeomorphic. Then 
$$\mc A(\partial M) = \mc A(\partial N) \ .$$ 
\end{thm} 

This is an immediate consequence of the  result below, where the Floer
homology groups are defined according to \cite{CFHW}. Indeed, this
definition of the Floer homology groups for a manifold $M$ 
depends only on the interior $\dot{M}$ 
(the admissible Hamiltonians have
compact support in $\dot{M}$) and on the possibility to ``complete''
$\dot{M}$ by adding a boundary of contact type.  

\begin{thm} Let $M$ be a manifold satisfying the  hypotheses of
  Theorem \ref{thm stab} and
  let $a\in \RR^*$. 

  a) The
  groups $FH^*_{]a-\epsilon, \, a+\epsilon]}(M)$ become independent of
  $\epsilon>0$ as soon as the latter is sufficiently small. We denote
  them by $FH^*_a(M)$. 

  b) If $(a, \, k) \notin \mc A(\partial M)$ for all $k\in \ZZ$ then
  $FH^k_a(M)=0$. If $(a, \, k) \in \mc A(\partial M)$ has multiplicity
  $n$, then it gives rise to $n$ copies of $\ZZ$ as direct summands in
  $FH^k_a(M)$ and to $n$ copies of $\ZZ$ as direct summands in
  $FH^{k+1}_a(M)$.  
\end{thm}

\medskip 

We note that the hypothesis $a\neq 0$ is necessary for the
groups $FH^*_{]a-\epsilon, \, a+\epsilon]}(M)$ to be well defined (cf. 
\S\ref{CiFH}).

\medskip 

{\it \small Sketch of proof.} The nondegeneracy hypothesis implies
that the characteristics have geometrically isolated images. In
particular, the values of their action form a discrete set $\mc A
\subset \RR$. If $a\in \mc A$ then it is possible to choose $\epsilon
>0$ such that $]a-\epsilon, \, a+\epsilon] \cap \mc A = \{ a \}$ and
the invariance property of symplectic homology ensures that all
$FH^*_{]a-\epsilon', \, a+\epsilon']}(M) $ are isomorphic for
$0<\epsilon' \le \epsilon$. The same holds if $a\notin \mc A$, with
the additional information that this group is now equal to zero by
definition. This proves a) and the first assertion in b). 

Let us now focus on the second assertion in b). The Liouville flow
gives rise to the usual trivialization of a neighbourhood of the
boundary as $\partial M \times [1-\delta, \, 1]$, $\delta >0$
sufficiently small. Denote $S$ 
the coordinate on the second factor. 
One considers  a cofinal family of Hamiltonians (Figure \ref{cofi
  spectrum}) which verify $H(p,\,
S) = h(S)$ on $\partial M \times [1-\delta, \, 1]$, with $h \le 0$
increasing, $h\equiv \tx{ct.}$ on $[1-\delta, \, r_1]$, $h$
strictly convex on $[r_1, \, r_2]$, $h' \equiv \tx{ct.}$
on $[r_2, \, r_3]$, $h$ strictly concave on $[r_3, \, r_4]$, $h\equiv
0$ on $[r_4, \, 1]$. Here $1-\delta < r_1 < r_2 < r_3 < r_4 <1$ and
$r_1 \longrightarrow 1$ as $H \stackrel \prec \longrightarrow
-\infty$.  
Moreover, one
assumes that the value of $h'$ on $[r_2, \, r_3]$ is not equal to the
action of any closed characteristic on $\partial M$.

\begin{figure}[h]  
      \begin{center}   
      \includegraphics{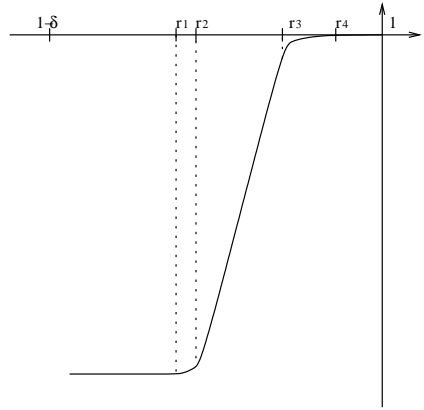}  
      \end{center}   
      \caption{Cofinal family in the proof of the stability of the
        action spectrum\label{cofi spectrum}}   
\end{figure}   

These Hamiltonians have highly degenerate $1$-periodic orbits but
the ones having nonzero action 
are made nondegenerate by an arbitrarily small perturbation
away from $\partial M \times [r_4, \, 1]$. 
After perturbing, 
the $1$-periodic orbits of $H$ fall into the following
classes: 
\begin{itemize} 
  \item constants in $\partial M \times [r_4, \, 1]$ with action $0$;
  \item isolated constants away from $\partial M \times [1-\delta, \,
  1]$, with action approximately equal to $-\min H \longrightarrow +\infty$;
 \item every characteristic $x$ with action less than $\max h'$ is seen
 twice in the unperturbed Hamiltonian: once in $\partial M \times
 [r_1, \, r_2]$ (call the corresponding orbit $x_1$) and once in
 $\partial M \times [r_3, \, r_4]$ (call the corresponding orbit
 $x_2$). After perturbation, each of these produces two nondegenerate
 periodic orbits $x_1^\pm$, $x_2^\pm$.  Moreover $A_H(x_1^\pm)
 \longrightarrow +\infty$, $A_H(x_2^\pm) \longrightarrow A(x)$. 
\end{itemize}

A computation of the Robbin-Salamon indices very similar in spirit to
the one already performed in \S\ref{coh Floer boule} shows that 
$$i_{RS}(x_1) = i_{RS}(x)+\frac 1 2 , \qquad  i_{RS}(x_2) =
i_{RS}(x)-\frac 1 2 \ ,$$
while 
$$i_{RS}(x_i^\pm)= i_{RS}(x_i)\mp \frac 1 2, \qquad i=1,\, 2 \ .$$

In the limit the only orbits whose action belongs to $[A(x) -\epsilon,
A(x)+\epsilon]$ are $x_2^\pm$, with indices $i_{RS}(x_2^-)=i_{RS}(x)$, 
$i_{RS}(x_2^+)=i_{RS}(x)-1 $ or, otherwise stated

$$-i_{RS}(x_2^-)=-i_{RS}(x), \qquad - i_{RS}(x_2^+)=-i_{RS}(x)+1 \ .$$

If, for a given value $a$ of the action, there is a single
characteristic taking on this value, then the computation of
$FH^*_{]a-\epsilon, \, a+\epsilon]}(H)$ is reduced to the
understanding of the differential in the Floer complex $\ZZ\langle
x_2^- \rangle \longrightarrow \ZZ \langle x_2^+ \rangle$. We have
seen in \S\ref{coh Floer boule} as a by-product of the direct
calculation of the Floer cohomology that the differential between
the perturbed images of a sphere of periodic orbits is zero. A direct
argument based on the implicit function theorem and the isotopy
invariance of Floer cohomology 
can be applied in the present situation to yield the same
result, and generalize it to the case where there are several closed
characteristics with action $a$. 

\hfill{$\square$}

\subsection{Applications to Weinstein's conjecture and obstructions to
  exact Lagrange embeddings.}  \label{Weinstein}
In this section we use the Floer cohomology groups as defined by
Viterbo (cf. section \ref{Homologie symplectique Viterbo}). 
We have already mentioned 
that they are useful for proving the existence
of closed characteristics on contact type hypersurfaces. The first key
tool is provided by the map 
\begin{equation} \label{morphisme Floer vers sing} 
\xymatrix{ FH^*(M) \ar[r]^{c^* \quad } &  H^{n+*}(M, \, \partial M) },
\quad n  
= \frac 1 2 \dim M  \ , 
\end{equation}  
obtained by restricting the range of the action and by using the isomorphism 
$$FH^*_{(a, \, 0^+]}(M) \simeq H^{n+*}(M, \, \partial M), \quad a<0 \ .$$  
Here  $0^+$ stands for a small enough positive number. 
As we have already pointed out in \S\ref{hom Floer var fermees}, when
restricting the range of the action and considering autonomous
Hamiltonians $H$ that are sufficiently small in the $C^2$ norm, the
Floer complex reduces to the Morse complex of the vector field 
$\nabla H$. By construction, the latter is outward pointing along
$\partial M$. We get in this way the homology of $M$ 
relative to the boundary, graded by the Morse index of the critical
points of $-H$, the function for which $\nabla H$ is negative
pseudo-gradient (see also \cite{teza mea}, Ch. 3 for further
details on this point). 

By the very construction of Floer homology we see that the
failure of $c^*$ to be bijective implies the existence of a closed
characteristic on  $\partial M$, which corresponds to the appearance
in the Floer complex of a generator other than a critical point of a
Morse function defined on 
$\tx{int}(M)$. We follow  
\cite{functors1} and consider the following definition.    
\begin{defi} {\sf \cite{functors1}} \label{formes algebriques conj Weinstein}  
A manifold $M$ satisfies the algebraic Weinstein conjecture 
 (AWC) if there is a ring of coefficients such that the map 
 $\xymatrix{ FH^*(M) 
\ar[r]^{c^* \quad } &  H^{n+*}(M, \, \partial M) }$ is not an isomorphism. 
\end{defi}

\medskip 

It is clear that  AWC implies the validity of Weinstein's conjecture for  
$\partial M$. In certain situations that we are going to state below,
the AWC property is inherited by codimension $0$ submanifolds in $M$
and it will prove useful to distinguish between the situations where
the morphism  $FH^*(M) \longrightarrow H^{n+*}(M, \, \partial M)$
fails to be injective or fails to be 
surjective (we shall then say that  $M$ satisfies case
a), respectively b) of the  AWC property). 

The second key tool is the existence 
\cite{functors1} \S2 of a transfer morphism   
\begin{equation} \label{morphisme transfer}  
\xymatrix{ FH^*(W) \ar[r]^{Fj^!} & FH^*(M) } \ , 
\end{equation}  
that is associated to a codimension $0$ inclusion 
$W \stackrel{j}  
\hookrightarrow M$, where $M$ has contact type boundary and  $\partial
W$ is of restricted contact type {\it in $W$}, with the meaning that
the symplectic form admits a primitive defined on the whole of  $W$ 
such that the Liouville vector field is transverse to 
$\partial W$. 

There are two remarks to be made about this transfer morphism. The
 first remark is that $Fj^!$ is defined only if the manifold 
$W$ verifies an additional
condition on the Floer trajectories (\cite{functors1}, 
p. 1000). Following Viterbo, we shall refer to this condition as {\it
 condition (A)}. It states that 
Floer trajectories {\it in $M$} for Hamiltonians of the type
 below, running between nonconstant $1$-periodic
orbits contained 
in $W$, must be {\it entirely contained
in $W$}. The main
point about this condition is that if it is violated
then there obviously is a periodic orbit on $\partial W$.

The second remark is that, if $\partial W$ carries no closed
characteristic, then the hypothesis of being of restricted contact
type in $W$ can be relaxed to contact type, and the morphism 
$Fj^!$ is still defined (condition (A) being assumed).

Here is a brief description of the morphism 
(\ref{morphisme transfer}) (see also Figure \ref{figure morphisme
  transfer}).  
A neighbourhood of  $\partial W$ in $W$ is trivialised by the Liouville
flow of  $W$ as  
$\partial W \times ]0, \, 1]$, due to the {\it restricted} contact
type nature of $\partial W$. Let us denote by $S_{_W}$ the
coordinate on the second factor. A neighbourhood of  $\partial M$ is
trivialised in the usual way by the Liouville flow of $M$ as 
$\partial M \times ]1-\delta, \, 1]$, $\delta >0$ 
and we denote by $S$ the coordinate on the second factor. The morphism 
 $Fj^!$ is constructed with the help of a family  $H=H_{\mu, \,
\lambda, \, \epsilon}$  of Hamiltonians that take into account the
characteristics on  $\partial W$ {\it as well as those on} $\partial
M$ and which bear the following form: $H=h_\lambda(S_{_W})$, $\epsilon
\le S_{_W}  \le 1$,  $h'_\lambda =\lambda$ and $H=k_\mu(S)$, $S\ge 1$,
$k'_\mu=\mu$, while $H$ is $C^2$-close to $0$ on $W \setminus \partial W \times
[\epsilon, \, 1]$ and $H=\lambda(1-\epsilon)$ on $M\setminus W$. The
transfer morphism is induced by the truncation morphism (for a
suitable almost complex structure $J$)
$$FH^*_{]0^-, \, b]}(H, \, J) \longrightarrow FH^*_{]a, \, b]}(H, \,
J) \ ,$$
with $a=a(\epsilon, \, \lambda)<0$ negative enough and 
$b=b(\epsilon, \, \lambda)>0$ positive enough. An isotopy 
argument identifies the (inverse) limit following $\epsilon$,
$\lambda$ and $\mu$ in the second term with  $FH^*(M)$. The hypothesis
mentioned above allow for a suitable choice of the
parameters $\epsilon$ and $\lambda$ for which the only orbits with
positive or close to zero action are the ones in $W$, and 
the first term is
identified with  $FH^*_{]0^-, \, b/\epsilon]}(K_\lambda)$,
$K_\lambda = k_\lambda(S_{_W})$, $S_{_W}\ge 1$ in $\widehat{W}$,
$K_\lambda=0$ on $W$. By taking the  (inverse) limit following 
$\lambda$ and $\epsilon$  one gets $FH^*(W)$. 

  \begin{figure}[h]  
        \begin{center}   
        \includegraphics{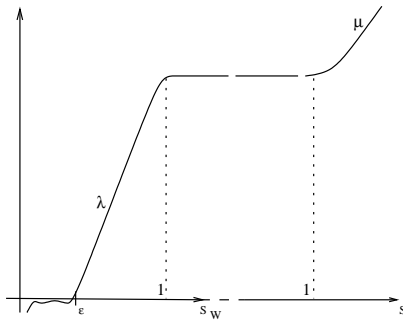}  
        \end{center}   
        \caption{Construction of the transfer morphism \cite{functors1} 
          \label{figure morphisme transfer}}    
  \end{figure}

The fundamental property of the transfer morphism is the commutativity
of the following diagram, which is proved by a careful analysis of the
various truncation morphisms which are involved
 \cite{functors1} Thm. 3.1. :  
\begin{equation} \label{diagramme commutatif morphisme transfer}   
\xymatrix{ FH^*(W) \ar[r]^{Fj^!} \ar[d]_{c^*} & FH^*(M) \ar[d]^{c^*}  
\\  
H^{n+*}(W, \, \partial W) \ar[r]^{j^!} & H^{n+*}(M, \, \partial M)   
} 
\end{equation}  
This allows one to deduce the following:
\begin{thm} {\sf (\cite{functors1}, Thm. 4.1) (Weinstein conjecture)} 
\label{application 
    Weinstein} Let $M$ be a connected symplectic manifold of dimension
 $2n$ such that the morphism  
$\xymatrix{ FH^{n} (M) 
\ar[r]^{c^* \quad} & H^{2n}(M, \, \partial M)}$ is not surjective. Any  
submanifold  $W \stackrel{j} \hookrightarrow M$ with 
 contact type boundary admits a closed characteristic on 
 $\partial W$. 
\end{thm}  
\demo Assume that $\partial W$ carries no closed characteristic. Then
condition (A) is trivially satisfied and the transfer morphism is
defined under the relaxed contact type assumption. 
The commutative diagram 
 (\ref{diagramme commutatif morphisme 
transfer}) therefore holds true. 
The morphism  $H^{2n}(W, \, \partial W) \longrightarrow 
H^{2n}(M, \, \partial M)$ is bijective if  $W$ is connected or
surjective if   
$W$ is disconnected. The hypothesis implies therefore that 
$FH^{n} (W) \longrightarrow H^{2n}(W, \, \partial W)$ is not
surjective. We infer the existence of a closed characteristic on 
 $\partial W$ and this gives a contradiction and proves that there is
always a closed characteristic on $\partial W$.   

\hfill{$\square$} 

\medskip 

 {\bf Remark.} If one uses field coefficients and the
manifold $M$ is connected, the nonsurjectivity of the 
morphism $FH^{n}(M) \longrightarrow H^{2n}(M, \, \partial M)$
is equivalent to its vanishing. 

\medskip 
 
Let us now mention some interesting cases of manifolds which satisfy
case b) of the AWC. 

\medskip 

1. We have already seen that
$$FH^*(D^{2n}) =0, \quad D^{2n}=\{z\in \CC^n \ : \ |z| \le 1 \} \ .$$ 
Through the preceding result, this directly implies Weinstein's
conjecture for hypersurfaces of contact type in 
 $\CC^n$ ( the original proof appeared in \cite{IHP}).  

2. It is proved in  \cite{functors1}, Thm. 4.2. 
that a subcritical Stein
manifold of dimension  $2n$ 
(i.e. having the homotopy type of a 
CW-complex of dimension  $n-1$) always satisfies case  b) of the  AWC
property: the morphism  $FH^{n}(M) \stackrel {c^*}
\longrightarrow H^{2n}(M, \, \partial M)$ is not surjective. The same
argument as above proves the existence of a closed characteristic on
any hypersurface of  contact type in a subcritical Stein manifold. 

3. The above result can also be deduced as 
   a consequence of a K\"unneth formula valid
   in Floer cohomology with field coefficients for a product of
   manifolds with restricted contact type boundary \cite{teza mea},
   Ch.2. There is a 
   commutative diagram where the horizontal arrows are
   isomorphisms:

\begin{equation} \label{Kunneth cohomologie}   
   {\scriptsize  
\xymatrix  
@C=13pt  
@R=30pt@W=1pt@H=1pt  
{ \bigoplus_{r+s = k} FH^r(M) \otimes    
  FH^s(N) \ar[rr]^{\quad \quad \ \ \sim} \ar[d]_{c^*\otimes c^*}   
& & FH^k(M \times N) \ar[d]^{c^*}  \\ 
\bigoplus_{r+s = k} H^{m+r}(M, \, \partial M) \otimes    
  H^{n+s}(N, \, \partial N)  \ar[rr]^{\quad \quad \ \ \sim} & &  
H^{m+n+k}(M \times N, \, \partial (M \times  
N) )  
}     
}  
\end{equation}

 Now a theorem of K. Cieliebak \cite{C} ensures that a subcritical Stein
 manifold $M$ is deformation equivalent to a split manifold $(N\times \CC,
 \, \omega_N \oplus \omega_{\text{\rm std}})$. The fact that
 $FH^*(\CC)=0$ allows one to infer the stronger result $FH^*(M)=0$. 

The commutativity of the diagram (\ref{Kunneth
  cohomologie}) has not been used in the above argument. Nevertheless
  it plays a crucial role in the proof of the stability of AWCb) 
  under products in
  the class of restricted contact type symplectic manifolds. 

4. I know of only one more case of manifold satisfying case b) of AWC,
   namely unit disc bundles associated to hermitian
   line bundles $\mc L$ with negative Chern class
   $c_1(\mc L) = -\lambda [\omega]$, $\lambda > 0$ over a symplectic
   base $(B, \, \omega)$. These are endowed with a symplectic form that
   restricts to the area form in the fibers and which equals the
   pull-back of $\omega$ on the horizontal distribution of a hermitian
   connection. One can compute in this case $FH^*(\mc{L})=0$ as a
   consequence of the existence of a spectral sequence valid in Floer
   homology \cite{teza mea}, Ch. 4.

\medskip 

Let us mention at this point that 
the only other available explicit computation of Floer homology concerns
cotangent bundles 
(\cite{coh cotangent, Sal coh cotg}) :  
$$FH^*(DT^*N) \simeq H^{*}(\Lambda N) \ ,$$ 
where  $N$ is a closed Riemannian manifold, 
$DT^*N =  \{ v\in T^*N \ : \ 
|v|\le 1 \}$ and $\Lambda N$ is the loop space of  $N$.  
The morphism $c^*$ coincides, modulo  Thom's isomorphism, with
the {\it surjection} $H^{*}(\Lambda N) \longrightarrow H^{*}(N)$ induced by
the inclusion  $N \hookrightarrow \Lambda N$. In this situation the
argument of Theorem \ref{application Weinstein} does not work
anymore in order to prove the Weinstein conjecture in cotangent
bundles, but one can nevertheless use the fact that the morphism
$FH^{n}(DT^{*}N) \longrightarrow H^{2n}(DT^{*}N, \, \partial DT^{*}N)
\simeq H^{n}(N)$ 
{\it is}
surjective in order to prove that there are no exact Lagrangian
embeddings $L\hookrightarrow M$ for manifolds $M$ satisfying AWCb) in
maximal degree. Here, by exact Lagrangian embedding one means that the
Liouville form $\lambda$ on $M$ restricts to an exact form on $L$. 
Indeed, according to a theorem of Weinstein, an exact Lagrangian
embedding would yield an embedding $DT^{*}L \stackrel j
{\hookrightarrow} M$ with $\partial DT^{*}L $ of restricted contact
type {\it in $M$}. One can prove that condition (A) is {\it always}
  satisfied  in this situation and therefore one would get a
  commutative diagram 
 $$\xymatrix{FH^n(DT^*L) \ar[r]^{Fj^!} \ar[d]_{c^*} & FH^n(M)
 \ar[d]^{c^*} \\
   H^{2n}(DT^*L, \, \partial DT^*L) \ar[r]^{\quad j^!} & H^{2n}(M, \,
 \partial M) 
}
 $$
This would imply that $FH^n(DT^*L) \longrightarrow H^{2n}(DT^{*}N, \,
\partial DT^{*}N)$ is not surjective, a contradiction. We have thus
proved the following theorem. 

\begin{thm} Let $M$ be a manifold with restricted contact type
  boundary which verifies AWCb) in maximal degree. There is no exact
  Lagrange embedding $L \hookrightarrow M$.
\end{thm}

\medskip 

The above theorem applies in particular for $M$ a subcritical Stein
manifold. We also note that the Weinstein conjecture in cotangent bundles is
proved in \cite{functors1} for simply connected manifolds with the
help of an equivariant version of Floer homology.

\medskip 
 
Before closing this section, let us remark that the preceding results
have their analogues in a homological setting. There still is  a
morphism 
\begin{equation} \label{morphisme sing vers Floer homologie} 
\xymatrix{ H_{n+*}(M, \, \partial M) \ar[r]^{\quad c_*} &  FH_*(M)} 
\end{equation} 
obtained by restricting the range of the action, as well as a transfer
morphism 
$$\xymatrix{ FH_*(M) \ar[r]^{Fj_!} & FH_*(W) } $$ 
which is defined under the same hypothesis as above. This fits into
the commutative diagram 
\begin{equation} \label{diagramme commutatif morphisme transfer 
    homologie}    
\xymatrix{ FH_*(W)  & FH_*(M) \ar[l]_{Fj_!} \\  
H_{n+*}(W, \, \partial W) \ar[u]^{c_*} & H_{n+*}(M,\, 
\partial M) \ar[u]_{c_*}   \ar[l]_{j_!} 
} 
\end{equation} 
The proof of these claims is dual to the one in
\cite{functors1}.

\section{Further reading and conclusions}

We have described three constructions and gave an application for each
of them. The underlying idea - we have repeated it over and over - is
to catch symplectically invariant information about the
characteristics in terms of Hamiltonians defined on the whole
manifold. A first theme to keep in mind is that there is no
``best''construction and that each of the ones we have presented has
virtues and  shortcomings, which arise from the different
behaviours at infinity or near the boundary that one imposes on the
admissible Hamiltonians. 

\medskip 

Let us mention some material for further reading - here ``further''
does not mean ``harder'', but only ``different and exploiting similar
ideas''. 

The paper of P. Biran, L. Polterovich and D. Salamon
\cite{BPS} develops at least three directions that are closely related
to the ideas that we have discussed. The first of them is the
construction of Floer homology groups based on the $1$-periodic orbits
belonging to a given free homotopy class of loops which is not
necessarily trivial. The second direction concerns the modification of
the class of admissible Hamiltonians. The authors prove existence
results for $1$-periodic orbits of compactly supported Hamiltonians in 
cotangent bundles of flat tori and in cotangent bundles of negatively
curved manifolds, under the only assumption that their value is
prescribed and sufficiently large on the zero section. There is an 
interplay between direct and inverse limits for this class of
admissible Hamiltonians and the class of compactly 
supported ones, which
allows the authors to construct {\it relative symplectic
  capacities} for the above types of cotangent bundles. Their
computation uses the results of Pozniak \cite{Po} on 
Floer homology for Hamiltonians
whose periodic orbits appear in Morse-Bott nondegenerate families,
like the ones in section \ref{coh Floer boule} above. 

A generalization of the above results for arbitrary cotangent bundles
over closed manifolds is in progress by J. Weber \cite{W}, using the
computation of truncated 
Floer homology groups for cotangent bundles \cite{coh
  cotangent, Sal coh cotg}. 

The paper by U. Frauenfelder and F. Schlenk \cite{FrSc} studies
the dynamics of compactly supported Hamiltonians  
on completions of compact manifolds with contact
type boundary or products of such objects - {\it split-convex
  manifolds}. The main tool is a kind of symplectic capacity, or
selector, associating to any Hamiltonian a value in its spectrum. The
selector is defined as the smallest value for which the image of the
fundamental class of the manifold through
a truncation morphism in Floer homology 
vanishes.  One interesting feature of the authors' construction is that
the emphasis is put on the Hamiltonians rather than on some geometric
levels, so that they need to consider only slow growths at
infinity.  Applications include existence of $1$-periodic
orbits for Hamiltonians having a displaceable support. The main point
is that displaceability implies uniform boundedness of the selector
for all the iterates of the corresponding Hamiltonian diffeomorphism. 

A recent survey of the Weinstein conjecture, including details on the
selector method, is provided by the paper of V.L. Ginzburg \cite{Gi}. 

\medskip 

Let us mention two directions for future research that seem
promising. The first concerns the class of admissible 
Hamiltonians prescribed for each of the homology theories that we have
presented. Their behaviour near the boundary
\cite{CFH, CFHW} or at infinity \cite{FH, functors1} 
is used in a crucial manner in order to obtain a priori bounds on the
Floer trajectories but, in some sense, is too rigid. This becomes
appearant as soon as one tries to perform geometric constructions on
Floer (co)homology, related to additional geometric structure on the
underlying manifold. As an example, the difficulty of the proof of the
K\"unneth formula (\ref{Kunneth cohomologie}) lies in the fact that
the componentwise sum
of two Hamiltonians that are linear at infinity is no longer linear at
infinity on the product manifold. It seems
to me of interest to enlarge the class of admissible
Hamiltonians in the setting of \cite{functors1}, 
the main point  being to still be able to prove a
priori $C^0$ bounds. An extension in this sense was accomplished in
the dissertation \cite{teza mea} through Hamiltonians that are
asymptotically linear at infinity. This points to the
investigation of new instances where the maximum principle can be
applied, for example by allowing the almost complex structure 
to vary as well at
infinity. 
 
A second interesting direction of investigation is to clarify the
relationship between (truncated) Floer homology and contact
homology of the boundary. Contact homology is a recent and powerful
invariant that can be defined intrinsically for any {\it contact
  manifold}, in particular for the boundary of contact type of a
symplectic manifold. This could lead to obstructions on the topology
of symplectic fillings of contact manifolds.

\bigskip 

\bigskip

{\sc Alexandru  Oancea, Department of Mathematics, ETHZ, 8092
  Z\"urich, Switzerland. } 
{\it E-mail address}: {\tt oancea\,@\,math.ethz.ch}

\end{document}